\documentclass[11pt]{article}
\usepackage{amssymb}

\newtheorem{definition}{Definition}[section]
\newtheorem{theorem}[definition]{Theorem}
\newtheorem{lemma}[definition]{Lemma}
\newtheorem{corollary}[definition]{Corollary}
\newtheorem{remark}[definition]{Remark}
\newtheorem{example}[definition]{Example}
\newtheorem{conjecture}[definition]{Conjecture}
\newtheorem{problem}[definition]{Problem}
\newtheorem{note}[definition]{Note}

\def\alg{\mathcal A}

\def\K{\mathbb K}

\newcommand{\beast}{\begin{eqnarray*}}
\newcommand{\eeast}{\end{eqnarray*}}

\begin{document}
\newenvironment{proof}{\noindent{\it Proof\/}:}{\par\noindent $\Box$\par}

\title{ \bf Two linear transformations each 
tri-
\\ 
\bf diagonal  
  with respect to an eigenbasis \\
\bf of the other; the $TD$-$D$ canonical form\\
 and the $LB$-$UB$ canonical form}
\author{Paul Terwilliger}
\date{}
%to get date printout, comment out above line 
\maketitle

\begin{abstract} 
Let $\K$ denote a field 
 and let $V$ denote a  
vector space over $\K$ with finite positive dimension.
We consider an ordered pair of linear
  transformations 
$A:V\rightarrow V$ and $B:V\rightarrow V$ 
which  
 satisfy both (i), (ii) below. 
\begin{enumerate}
\item There exists a basis for $V$ with respect to which
the matrix representing $A$ is irreducible tridiagonal and
the matrix representing $B$ is diagonal.
\item There exists a basis for $V$ with respect to which
the matrix 
representing $A$ is diagonal and 
the matrix representing $B$ is irreducible tridiagonal.
\end{enumerate}
We call such a pair a {\it Leonard pair} on $V$.
We introduce two canonical forms for Leonard
pairs. We call these the
 {\it $TD$-$D$ canonical form} and the
{\it $LB$-$UB$ canonical form}.
In the 
 $TD$-$D$ canonical form 
the Leonard pair 
is represented by an irreducible tridiagonal matrix
and a diagonal matrix, subject to a certain normalization.
In the 
$LB$-$UB$ canonical form
the  Leonard pair 
is represented by a lower bidiagonal matrix and an upper bidiagonal matrix,
subject to a certain normalization.
We describe the two canonical forms in
 detail. 
As an application we obtain the following results.
Given square matrices $A,B$ over $\K$,
with $A$ tridiagonal and $B$ diagonal,
we display a necessary and  sufficient condition for
$A,B$ to represent  a Leonard pair.
Given square matrices $A,B$ over $\K$, with
$A$ lower bidiagonal and $B$ upper
bidiagonal,
we display a necessary and  sufficient condition for
$A,B$ to represent  a Leonard pair.
We briefly discuss how Leonard pairs correspond to 
the $q$-Racah polynomials
and some related polynomials in the Askey scheme.
We present some open problems concerning
Leonard pairs.
\end{abstract}

\section{Introduction}

\medskip
\noindent 
We begin by recalling the notion of a {\it Leonard pair}
\cite{TD00}, \cite{Terint}, \cite{rome}, 
\cite{LS24}, \cite{LS99}, \cite{qSerre}.
We will use the following terms.
Throughout this paper, 
when we refer to a matrix, we mean a square matrix.    
 A matrix is called {\it tridiagonal} whenever
  each nonzero entry lies on either the
diagonal, the subdiagonal, or the superdiagonal.
A 
tridiagonal matrix is called {\it irreducible} 
whenever each entry on the subdiagonal is nonzero and 
 each entry on the superdiagonal is nonzero.

\medskip
\noindent
We now define a Leonard pair.
For the rest of this paper 
$\K$ will denote a field.

\begin{definition}
 \cite{LS99}
\label{def:lprecall}
\rm 
Let  
 $V$ denote a  
vector space over $\K$ with finite positive dimension.
By a {\it Leonard pair} on $V$
we mean an ordered pair of linear transformations 
$A:V\rightarrow V$ and $A^*:V\rightarrow V$ 
which  
satisfy both (i), (ii) below.
\begin{enumerate}
\item There exists a basis for $V$ with respect to which
the matrix representing $A$ is irreducible tridiagonal and
the matrix representing $A^*$ is diagonal.
\item There exists a basis for $V$ with respect to which
the matrix 
representing $A$ is diagonal and 
the matrix representing $A^*$ is irreducible tridiagonal.
\end{enumerate}
\end{definition}
\noindent 

\begin{note} 
\rm 
According to a common notational convention
$A^*$ denotes the 
conjugate-transpose of
 $A$. We emphasize we are not using
this convention. In a Leonard pair $A,A^*,$ the linear
transformations $A$ and $A^*$ are arbitrary subject to
(i), (ii) above.
\end{note}

%
%Our use of the name ``Leonard pair'' is motivated by a connection
%to a theorem of D. Leonard 
%\cite[p.  260]{BanIto}, \cite{Leon}
%involving the $q$-Racah and related
%polynomials of the Askey Scheme. See 
%\cite[App. A]{LS99} and 
%\cite{LS24} for more information
%on this.

\medskip
\noindent We give some background
on Leonard pairs. There is a  correspondence between Leonard
pairs and a family of orthogonal polynomials  
consisting of the $q$-Racah polynomials and some related polynomials
in the Askey scheme.
This correspondence is discussed 
in 
\cite{rome},
\cite{LS24},
\cite[Appendix A]{LS99} 
and in Section
27 below.
%This correspondIt is also implicit in
%\cite[p.  260]{BanIto}, \cite{Leon}.
The reference \cite{KoeSwa} contains detailed information about
the Askey scheme.

\medskip
\noindent Leonard pairs play a role in representation theory.
For instance Leonard pairs arise naturally
in the representation theory 
 of the Lie algebra $sl_2$ 
 \cite{TD00}, the quantum algebra
$U_q(sl_2)$ 
\cite{Koelink3},
\cite{Koelink1},
\cite{Koelink2}, 
\cite{Koelink4},
\cite{koo3},
\cite[Chapter 4]{Hjal},
\cite{Terint}, 
\cite{qSerre},
the Askey-Wilson algebra 
\cite{GYZnature},
 \cite{GYLZmut},
%\cite{GYZTwisted},
\cite{GYZlinear},
%\cite{GYZspherical},
\cite{Zhidd} 
%\cite{ZheCart},
%\cite{Zhidden},
and 
the Tridiagonal algebra 
\cite{TD00},
\cite{LS99},
\cite{qSerre}.

\medskip
\noindent
Leonard pairs play a role  in combinatorics. For instance
Leonard pairs can be constructed from
certain  partially ordered sets
\cite{rome}.
Also, there exists a combinatorial object  
called a  $P$- and $Q$-polynomial
association scheme \cite{BanIto}, \cite{bcn}, \cite{Leopandq},
\cite{Tercharpq}, \cite{Ternew}.
Leonard pairs have been used to describe
certain irreducible modules for the subconstituent algebra of
 these association  schemes \cite{TersubI}, 
\cite{TersubII}, \cite{TersubIII}.  
See \cite{Cau}, \cite{CurNom}, \cite{Curspin}, \cite{go}, 
%\cite{HobIto},
\cite{TD00}
%\cite{Tan}
for more
information on Leonard pairs and association  schemes.

\medskip
\noindent 
Leonard pairs are closely related to the  work of 
Grunbaum and Haine on the 
 ``bispectral problem''
\cite{GH7},
\cite{GH6}.
See
\cite{GH4},
%\cite{GH5},
\cite{GH1}, 
\cite{GH3},
\cite{GH2} 
for related work. 

\medskip
\noindent
The rest of this introduction
contains a detailed summary of  the present paper.

\medskip
\noindent 
In this paper we introduce two canonical forms for
Leonard pairs.
The first of these is called
 the {\it $TD$-$D$ canonical form}.
In this form
the Leonard pair is represented
by an irreducible tridiagonal matrix and a diagonal matrix,
subject to a certain normalization.
To describe the second form we make a definition.
A matrix is said to be {\it lower bidiagonal} whenever
each nonzero entry lies 
on either the diagonal or the subdiagonal.
A matrix is said to be  {\it upper bidiagonal }
whenever its transpose is lower bidiagonal.
We call our second form the {\it $LB$-$UB$ canonical form}.
In this form
the Leonard pair is represented
by a lower bidiagonal matrix and an upper bidiagonal matrix, 
subject to a certain normalization.

\medskip
\noindent
We fix some notation. 
Let $d$ denote a nonnegative integer.
We let $\hbox{Mat}_{d+1}(\K)$ denote the
$\K$-algebra consisting of all
$d+1$ by $d+1$ matrices which have entries in $\K$. We
index the rows and columns by $0,1,\ldots, d$.
Any $\K$-algebra which is isomorphic to 
$\hbox{Mat}_{d+1}(\K)$ 
will be called a {\it  matrix algebra over $\K$ of diameter 
$d$}.

\medskip
\noindent Before proceeding we sharpen our concept of
a Leonard pair.
Let $\cal A$ denote a matrix algebra over $\K$ 
and let $V$ denote an irreducible
left $\cal A$-module.
By a {\it Leonard pair in $\cal A$} we mean 
an ordered pair of elements taken from $\cal A$
which act on $V$ as a Leonard pair in the sense of
 Definition 
\ref{def:lprecall}.
Let $A,A^*$ denote a Leonard pair in $\cal A$.
Then  $A$ and $A^*$ together generate  $\cal A$
\cite[Corollary 3.2]{LS24}.
By a {\it Leonard pair over $\K$} we mean
a sequence ${\cal A},  A,A^*$ where $\cal A$
is a matrix algebra over $\K$ and $A,A^*$ is
a Leonard pair in $\cal A$.
We call $\cal A$ the {\it ambient algebra}
and suppress it in the notation, referring only to $A,A^*$.
Let $A,A^*$ denote a Leonard pair over $\K$.
By the {\it diameter }
of this pair we mean the diameter of its ambient algebra.
By  the {\it underlying module} for this pair we mean an irreducible
left module for its ambient algebra.
For the rest of this Introduction, 
when we refer to a scalar we mean an element of $\K$.
When we refer to a Leonard pair it is assumed to be over $\K$.

\medskip
\noindent
We recall the notion of an {\it eigenvalue sequence} 
for a Leonard pair.
Let $A,A^*$ denote a Leonard pair.
By definition there exists a basis for
the underlying module
with respect to which the matrix representing  $A$ is diagonal
and the matrix representing $A^*$
is irreducible
tridiagonal. In the matrix representing $A$
the diagonal entries are the eigenvalues of $A$
and it turns out these are  
 mutually distinct
\cite[Lemma 1.3]{LS99}.
  Therefore the sequence of diagonal
entries gives an ordering of the eigenvalues of $A$.
We call this sequence an {\it eigenvalue sequence} for $A,A^*$.
Given an eigenvalue sequence for $A,A^*$, if we invert the
order of the sequence we get another
eigenvalue sequence for $A,A^*$.
Moreover $A,A^*$ has no  
further eigenvalue sequence.
To clarify this let $d$ denote the diameter of $A,A^*$.
Then $A,A^*$ has exactly two eigenvalue sequences if
$d\geq 1$ and
a single eigenvalue sequence if $d=0$.
By a 
 {\it dual eigenvalue sequence} for $A,A^*$ we mean an
 eigenvalue sequence for the Leonard  pair $A^*,A$.

\medskip
\noindent
A {\it Leonard system} is essentially a
Leonard pair, together 
with an eigenvalue sequence and a dual eigenvalue
sequence for that pair.
For the duration of this Introduction 
we take  this 
as the definition of a Leonard system.
(In the
main part of our paper we will define a Leonard system in a slightly
different manner
in which the eigenvalues are replaced by their corresponding
primitive idempotents.)

\medskip
\noindent 
We mentioned each Leonard system involves a Leonard pair;
we call this pair  the {\it associated} Leonard pair.
The set of Leonard systems associated with
a given Leonard pair will be called the {\it associate class} for
that pair. 
In order to describe the associate classes 
we use the following notation.
Let $\Phi$ denote a Leonard system.
If we invert the ordering
on the eigenvalue sequence of $\Phi$ we get a Leonard system 
which we denote by $\Phi^\Downarrow$.
If we instead invert the ordering on the dual eigenvalue
sequence of $\Phi$ we get a Leonard system which we denote 
by $\Phi^\downarrow$.
We view $\downarrow, \Downarrow $ as permutations on the set
of all Leonard systems. 
These permutations 
are commuting involutions
and therefore induce an action of the Klein 4-group
on the set of all 
 Leonard systems.
The orbits of this action coincide with the associate classes.

\medskip
\noindent We discuss the notion of {\it isomorphism} for Leonard
pairs and Leonard systems.
Let $A,A^*$ and $B,B^*$ denote Leonard pairs.
By an {\it isomorphism  of Leonard pairs from
 $A,A^*$ to $B,B^*$}
we mean an isomorphism of $\K$-algebras from the
ambient algebra of $A,A^*$ to the ambient algebra
of $B,B^*$ which sends $A$ to $B$  and $A^*$ to
$B^*$.
We say $A,A^*$ and $B,B^*$  are {\it isomorphic} whenever
there exists an isomorphism of Leonard pairs
from $A,A^*$ to $B,B^*$. We say
two given Leonard systems are {\it isomorphic} whenever
(i) their associated Leonard pairs are isomorphic;
(ii) their eigenvalue sequences coincide; and (iii) their  
dual  eigenvalue sequences coincide.

\medskip
\noindent The set of Leonard systems is partitioned into
both isomorphism classes and associate classes. 
These partitions are related as follows.
Let $A,A^*$ denote a Leonard pair and let $d$ denote
the diameter.
If $d\geq 1$ then 
the corresponding associate class contains
 four
Leonard systems and these are mutually nonisomorphic.
If 
$d= 0$ then the corresponding associate class contains a single
Leonard system.

\medskip
\noindent
Before proceeding with Leonard systems
we introduce the notion of a 
{\it parameter array}. A parameter
array is a finite sequence of scalars 
 which satisfy
a certain list of equations and inequalities. 
We care about parameter arrays because
it turns out they are in bijection with
the isomorphism classes of Leonard systems.
A parameter array is defined as follows.
Let $d$ denote a nonnegative integer.
By a {\it parameter array of diameter $d$} we mean a sequence
of scalars 
$(\theta_i, \theta^*_i, i=0..d;  \varphi_j, \phi_j, j=1..d)$
which satisfy (i)--(v) below.
\begin{enumerate}
\item $ \varphi_i \not=0, \qquad \phi_i\not=0 \qquad \qquad \qquad\qquad (1 \leq i \leq d)$.
\item $ \theta_i\not=\theta_j,\qquad  \theta^*_i\not=\theta^*_j\qquad $
if $\;\;i\not=j,\qquad \qquad \qquad (0 \leq i,j\leq d)$.
\item $ {\displaystyle{ \varphi_i = \phi_1 \sum_{h=0}^{i-1}
{{\theta_h-\theta_{d-h}}\over{\theta_0-\theta_d}} 
\;+\;(\theta^*_i-\theta^*_0)(\theta_{i-1}-\theta_d) \qquad \;\;(1 \leq i \leq d)}}$.
\item $ {\displaystyle{ \phi_i = \varphi_1 \sum_{h=0}^{i-1}
{{\theta_h-\theta_{d-h}}\over{\theta_0-\theta_d}} 
\;+\;(\theta^*_i-\theta^*_0)(\theta_{d-i+1}-\theta_0) \qquad (1 \leq i \leq d)}}$.
\item The expressions
\beast
{{\theta_{i-2}-\theta_{i+1}}\over {\theta_{i-1}-\theta_i}},\qquad \qquad  
 {{\theta^*_{i-2}-\theta^*_{i+1}}\over {\theta^*_{i-1}-\theta^*_i}} 
 \qquad  \qquad 
\eeast
are equal and independent of $i$ for $2 \leq i \leq d-1$.
\end{enumerate}

\medskip
\noindent 
We give  a bijection from the set of isomorphism
classes of Leonard systems to the set of parameter arrays.
Let $\Phi$ denote a Leonard system.
To  $\Phi$ we attach
the following four sequences of scalars.
The first two sequences are the eigenvalue sequence
of $\Phi$ and the dual eigenvalue sequence of $\Phi$.
Let us denote these by
$\theta_0, \theta_1,\ldots,\theta_d$ 
and $\theta^*_0, \theta^*_1,\ldots,\theta^*_d$, respectively.
By a slightly technical construction which we
omit for now, we obtain a third sequence of scalars
$\varphi_1, \varphi_2, \ldots, \varphi_d$.
We call this the {\it first split sequence}
of $\Phi$.
We let $\phi_1, \phi_2, \ldots, \phi_d$ denote
the first split sequence for $\Phi^{\Downarrow}$ and
call  this the {\it  second split sequence} of $\Phi$. 
By \cite[Theorem 1.9]{LS99} 
a sequence of scalars
$p=(\theta_i, \theta^*_i, i=0..d;  \varphi_j, \phi_j, j=1..d)$
is a parameter array if and only if there exists a Leonard system
$\Phi$ with eigenvalue sequence 
$\theta_0, \theta_1,\ldots,\theta_d$,
dual eigenvalue sequence 
$\theta^*_0, \theta^*_1,\ldots,\theta^*_d$, 
first split sequence 
$\varphi_1, \varphi_2, \ldots, \varphi_d$,
and second split sequence 
 $\phi_1, \phi_2, \ldots, \phi_d$. 
If $\Phi$ exists then $\Phi$ is unique up to isomorphism.
In this case  we call $p$ the {\it parameter array of $\Phi$}.
The map which sends a Leonard system to its parameter array
induces the desired bijection from the 
 set of isomorphism
classes of Leonard systems to the set of parameter arrays.

\medskip
\noindent Earlier we described an action of the Klein 4-group
on the set of Leonard systems.
The above bijection induces an action of the same group
on the set of parameter arrays. We now describe this action.
Let $\Phi$ denote a Leonard system and let
$p=(\theta_i, \theta^*_i, i=0..d;  \varphi_j, \phi_j, j=1..d)$
denote the parameter array of $\Phi$.
The parameter array of $\Phi^\downarrow $ is $p^\downarrow$
where
$p^\downarrow:=(\theta_i, \theta^*_{d-i}, i=0..d;  
\phi_{d-j+1}, \varphi_{d-j+1}, j=1..d)$.
The parameter array of $\Phi^\Downarrow$ is
$p^\Downarrow$ where
$p^\Downarrow:=
(\theta_{d-i}, \theta^*_i, i=0..d; 
\phi_{j}, \varphi_{j}, j=1..d)$ \cite[Theorem 1.11]{LS99}.

\medskip
\noindent
Let $A,A^*$ denote a Leonard pair.
By a {\it parameter array of $A,A^*$}
we mean the parameter array of an associated Leonard
system. We observe that if $p$ is a
parameter array of $A,A^*$ then so are 
$ p^\downarrow,  p^\Downarrow, p^{\downarrow \Downarrow}$
and $A,A^*$ has no further parameter arrays.
We comment on the distinctness of these arrays.
Let $d$ denote the diameter of $A,A^*$.
Then $ p, p^\downarrow,  p^\Downarrow, p^{\downarrow \Downarrow}$
are
mutually distinct if $d\geq 1$ and 
identical if $d=0$.
Therefore $A,A^*$ has exactly 
 four parameter arrays if $d\geq 1$  and just
one parameter array if $d=0$.

\medskip
\noindent 
We now describe    
the $TD$-$D$ canonical form.

\medskip
\noindent 
We define what it means for a given Leonard  system
to be in $TD$-$D$ canonical form.
Let $\Phi$ denote a Leonard system with
eigenvalue sequence  
 $\theta_0, \theta_1, \ldots, \theta_d$
and dual eigenvalue sequence 
 $\theta^*_0, \theta^*_1, \ldots, \theta^*_d$.
Let $A,A^*$
denote the associated Leonard pair.
Then $\Phi$ is in {\it   
$TD$-$D$ canonical form} whenever
(i) the ambient algebra of $A,A^*$ is 
$\hbox{Mat}_{d+1}(\K)$;
(ii) 
$A$ is tridiagonal and $A^*$ is diagonal;
(iii)
$A$ has constant row sum $\theta_0$ and $A^*_{00}=\theta^*_0$.

\medskip
\noindent 
We  describe the Leonard systems which are in
$TD$-$D$ canonical form.
In order to do this we consider
the set of parameter arrays.
We define two functions 
on this set.
We call these functions
$T$ and $D$.
Let
$p=(\theta_i, \theta^*_i, i=0..d;  \varphi_j, \phi_j, j=1..d)$
 denote a parameter array.
The image $p^T$ is the tridiagonal matrix in 
$\hbox{Mat}_{d+1}(\K)$
which has 
the following entries. The diagonal entries are
\beast
p^T_{ii} = \theta_i + \frac{\varphi_i}{\theta^*_i-\theta^*_{i-1}}
 + \frac{\varphi_{i+1}}{\theta^*_i-\theta^*_{i+1}}
\eeast
for $0 \leq i \leq d$, where $\varphi_0=0$, $\varphi_{d+1}=0$
and where $\theta^*_{-1}, \theta^*_{d+1}$ denote indeterminates.
The superdiagonal and subdiagonal entries are
\beast
p^{T}_{i-1,i}= 
\varphi_i \frac{\prod_{h=0}^{i-2}(\theta^*_{i-1}-\theta^*_h)
}
{\prod_{h=0}^{i-1}(\theta^*_{i}-\theta^*_h)
},
\qquad \qquad 
p^{T}_{i,i-1}=
 \phi_{i} \frac{\prod_{h=i+1}^d (\theta^*_i-\theta^*_h)
}
{\prod_{h=i}^d(\theta^*_{i-1}-\theta^*_h)
}
\eeast
for $1 \leq i \leq d$. 
The image 
$p^D$ is $\hbox{diag}(\theta^*_0, \theta^*_1, \ldots,\theta^*_d)$.
The significance of $T$ and $D$ is the following. 
Given a  Leonard system in $TD$-$D$ canonical form
the associated Leonard pair is $p^T,p^D$ where
$p$ denotes the corresponding parameter array.

\medskip
\noindent Let $\Phi$ denote a Leonard system. By a {\it $TD$-$D$
canonical form for $\Phi$}, we mean a Leonard system which is
isomorphic to $\Phi$ and which is in $TD$-$D$ canonical form.
 We show there exists a unique $TD$-$D$ canonical form for
 $\Phi$.

%may1kill
%\medskip
%\noindent 
%In order to clarify the situation so far
%we give a bijection from the set of Leonard systems which are
%in $TD$-$D$ canonical form, to the set of parameter arrays.
%The bijection sends each Leonard system to its parameter array.

\medskip
\noindent 
%Consider the set of Leonard systems which are in $TD$-$D$
%canonical form.
%We comment on how this set intersects
%the associate classes of Leonard systems. Each of these
%classes contains at most one Leonard system which is
%in 
% $TD$-$D$ canonical form.
%
%\medskip
%\noindent 
Let $A,A^*$ denote a Leonard pair and consider its
set of associated Leonard systems.
From the construction 
this set contains at most one Leonard system
which is in 
$TD$-$D$ canonical form.
The case in which this
Leonard system exists is of interest; to describe this case 
we define a $TD$-$D$ canonical form for Leonard pairs. 
We do this as follows.

\medskip
\noindent We define what it means for a Leonard pair to
be in $TD$-$D$ canonical form.
Let $A,A^*$ denote a Leonard pair and let
$\theta_0, \theta_1, \ldots, \theta_d$ denote an eigenvalue
sequence for this pair.
Then $A,A^*$ is 
in {\it $TD$-$D$ canonical form} whenever (i) the ambient
algebra of $A,A^*$ is 
$\hbox{Mat}_{d+1}(\K)$;
(ii) $A$ is tridiagonal and $A^*$ is diagonal;
(iii) $A$ has constant row sum and this sum is $\theta_0$ or
$\theta_d$.

\medskip
\noindent 
We just defined the $TD$-$D$ canonical form for Leonard
pairs and earlier we defined this form for 
Leonard systems.
These two versions
 are related as follows.
A given Leonard pair is in $TD$-$D$ canonical form if and only if
there exists an associated Leonard system which is in $TD$-$D$ 
canonical form.
%may1kill
%To state things from another point of view
%we give a bijection from the set of Leonard systems 
%which are in $TD$-$D$ canonical form, to the set of Leonard
%pairs which are in $TD$-$D$ canonical form.
%This bijection sends each Leonard system to the associated Leonard
%pair.
%
%\medskip
%\noindent 
%Consider the set of Leonard systems which are in $TD$-$D$ 
%canonical form.  
%So far we have given two bijections involving this set.
% We first gave a bijection from this set
%to the set  of parameter arrays.
%We then gave a bijection from
%this set 
%to the set of Leonard pairs which are in 
%$TD$-$D$ canonical form.
%Composing the inverse of the first bijection
%with the second, we get a bijection from the set of parameter arrays
%to the set of Leonard pairs which are in $TD$-$D$ canonical form.
%This bijection sends each parameter array $p$ to the pair $p^T,p^D$. 

\medskip
\noindent 
Let $A,A^*$ denote a Leonard pair. By a {\it $TD$-$D$
canonical form
for $A,A^*$} we mean a Leonard pair which is isomorphic to 
$A,A^*$ and which is in $TD$-$D$ canonical form.
We describe the $TD$-$D$ canonical forms for $A,A^*$.
To do this
we give a bijection from the set of parameter arrays for $A,A^*$
to the set of $TD$-$D$ canonical forms for $A,A^*$. This bijection
sends each parameter array $p$ to the pair $p^T,p^D$.
To clarify this let $d$ denote the diameter of $A,A^*$.
If $d\geq 1$ then
there exists exactly four  
 $TD$-$D$ canonical forms for $A,A^*$.
If $d=0$ then there exists a unique $TD$-$D$ canonical form for
$A,A^*$.

%\medskip
%\noindent 
%This takes us to the end of our formal argument. 

\medskip
\noindent 
We give several applications of our theory.
For the first application  
let $d$ denote a nonnegative integer and let
$A,A^*$ denote matrices in 
$\hbox{Mat}_{d+1}(\K)$.
We give a necessary and sufficient condition
for  $A,A^*$ to be a Leonard pair
in 
$\hbox{Mat}_{d+1}(\K)$ which is in $TD$-$D$ canonical form.
Indeed we show the following are equivalent:
(i)
the pair $A,A^*$ is a Leonard pair in
$\hbox{Mat}_{d+1}(\K)$ which is in
$TD$-$D$ canonical form; (ii)
there exists a parameter array $p$ of diameter $d$
such that $A=p^T$ and $A^*=p^D$.

\medskip
\noindent
Our second application is similar to the first but more general.
Again 
let $d$ denote a nonnegative integer and let
$A,A^*$ denote matrices in 
$\hbox{Mat}_{d+1}(\K)$.
Let us assume $A$ is tridiagonal  and $A^*$ is diagonal.
We give a necessary and sufficient condition for
 $A,A^*$ to be  
a Leonard pair in 
$\hbox{Mat}_{d+1}(\K)$.
This condition is given
in Theorem
\ref{thm:tdcrit}.

\medskip
\noindent This completes our description of the 
$TD$-$D$ canonical form. 
 Our description of the $LB$-$UB$ canonical form
runs along similar lines; we save the details for
the main body of the paper.
We comment that 
in the main body of the 
paper it will be convenient to treat the
$LB$-$UB$ canonical form before the
$TD$-$D$ canonical form.

\medskip
\noindent 
As we proceed through the paper we illustrate our results
using two running examples
which involve specific parameter arrays.

\medskip
\noindent Near the end of the paper
we  discuss how
Leonard pairs correspond
 to the $q$-Racah polynomials
and some 
related polynomials in the Askey scheme.
The general idea is the following.
Given a Leonard pair $A,A^*$ 
the corresponding polynomials 
give the entries in a transition matrix
which takes a basis satisfying
Definition
\ref{def:lprecall}(i)
to a basis satisfying
Definition
\ref{def:lprecall}(ii).
We compute these polynomials explicitly for our
two examples. 
For these examples the polynomials turn out to be Krawtchouk
polynomials and $q$-Racah polynomials.

\medskip
\noindent 
At the end of the paper we  present
some open problems concerning Leonard pairs.

\section{Leonard systems}
\medskip
\noindent We now begin our formal argument. 
Our first goal is to recall our working definition of a Leonard
system.
We begin with some notation.

\medskip
\noindent
Let $d$  denote  a  nonnegative
integer.
We let $\K^{d+1}$ denote the vector space over $\K$
consisting
of all $d+1$ by $1$  matrices which have entries in $\K$.
We index the rows by 
$0,1,\ldots, d$. We view $\K^{d+1}$ as 
a left
module for 
$\hbox{Mat}_{d+1}(\K)$ under matrix multiplication.
We observe this module is irreducible.
We let ${\mathcal A}$ denote a $\K$-algebra isomorphic to
$\hbox{Mat}_{d+1}(\K)$.
From now on when we refer to an ${\mathcal A}$-module we mean a left 
 ${\mathcal A}$-module. Let $V$ denote an irreducible $\cal A$-module.
We remark that $V$ is unique up to isomorphism of $\cal A$-modules,
and that $V$ has dimension $d+1$.
 Let $v_0, v_1, \ldots, v_d$ denote a basis for $V$.
For $X \in \mathcal A$ and 
 $Y \in 
\hbox{Mat}_{d+1}(\K)$,  we say $Y$ {\it represents $X$ with respect
to 
  $v_0, v_1, \ldots, v_d$} whenever 
$
X v_j =\sum_{i=0}^d Y_{ij}v_i
$
for $0 \leq j\leq d$.
For  $A\in \mathcal A$,
%By an {\it eigenvalue }
%of $A$ we mean a root of the minimal polynomial of $A$.
%The eigenvalues  of $A$ are contained in the algebraic closure of $\K$.
we say $A $  is {\it
multiplicity-free} whenever it has 
$d+1$ distinct  eigenvalues
in $\K$.
Assume $A$ is  multiplicity-free.
Let $\theta_0, \theta_1, \ldots, \theta_d$ denote an ordering of 
the eigenvalues
of $A$, and for $0 \leq i \leq d$   put 
\begin{equation}
E_i = \prod_{{0 \leq  j \leq d}\atop
{j\not=i}} {{A-\theta_j I}\over {\theta_i-\theta_j}},
\label{eq:primiddef}
\end{equation}
where $I$ denotes the identity of $\mathcal A$.
We observe
(i) $AE_i = \theta_iE_i$ $(0 \leq i \leq d)$;
%\label{eq:primid1S99}
(ii) $E_iE_j = \delta_{ij}E_i$ $(0 \leq i,j\leq d)$;
%\label{eq:primid2S99}
(iii) $\sum_{i=0}^d E_i = I$;
%\label{eq:primid3S99}
(iv) $A=\sum_{i=0}^d \theta_i E_i$.
Let $\mathcal D$ denote the subalgebra of $\mathcal A$ 
generated by $A$.
Using (i)--(iv) we find the sequence
$E_0, E_1, \ldots, E_d$ is a  basis for the $\K$-vector space 
$\mathcal D$.
We call $E_i$  the {\it primitive idempotent} of
$A$ associated with $\theta_i$.
It is helpful to think of these primitive idempotents as follows. 
Let $V$ denote an irreducible  $\mathcal A$-module. Then
\begin{eqnarray}
V = E_0V + E_1V + \cdots + E_dV \qquad \qquad (\hbox{direct sum}).
\label{eq:VdecompS99}
\end{eqnarray}
For $0\leq i \leq d$, $E_iV$ is the (one dimensional) eigenspace of
$A$ in $V$ associated with the 
eigenvalue $\theta_i$, 
and $E_i$ acts  on $V$ as the projection onto this eigenspace. 

\begin{definition}
\label{def:idseq}
\rm 
Let $d$ denote a nonnegative integer and let
$\mathcal A$ denote a $\K$-algebra isomorphic to 
$\hbox{Mat}_{d+1}(\K)$. 
Let $A,A^*$ denote an ordered pair consisting of
multiplicity-free  elements in $\mathcal A$.
By an {\it idempotent sequence} for $A,A^*$ we mean 
an ordering $E_0, E_1, \ldots, E_d$ of the primitive idempotents
of $A$ such that
\beast
E_iA^*E_j = \cases{0, &if $\;\vert i-j\vert > 1$;\cr
\not=0, &if $\;\vert i-j \vert = 1$\cr}
\qquad \qquad (0 \leq i,j\leq d).
\eeast
We observe that if 
 $E_0, E_1, \ldots, E_d$ is an
idempotent sequence for $A,A^*$   
then  so is 
 $E_d, E_{d-1}, \ldots, E_0$ and
  $A,A^*$ has no further idempotent sequence.
By a {\it dual idempotent sequence} for $A,A^*$ we mean
an idempotent sequence for $A^*,A$.
\end{definition}

\begin{definition} \cite{LS99}
\label{def:defls}
\rm 
Let $d$  denote  a  nonnegative
integer
 and let $\alg$ 
denote a $\K$-algebra isomorphic to 
$\hbox{Mat}_{d+1}(\K)$. 
By a {\it Leonard system} in $\alg$ we mean a 
sequence 
\begin{equation}
\; \Phi = (A, A^*; E_i, E^*_i,i=0..d )
\label{eq:ourstartingpt}
\end{equation}
 which satisfies  (i)--(iii) below.
\begin{enumerate}
\item Each of $A$,  $A^*$ is a multiplicity-free element of $\alg$.
\item $E_0,E_1,\ldots,E_d$ is an idempotent sequence for $A,A^*$.
\item $E^*_0,E^*_1,\ldots,E^*_d$ is a dual idempotent sequence for
$A, A^*$.
\end{enumerate}
We call $d$ the {\it diameter} of $\Phi$ and say 
$\Phi$ is {\it over } $\K$.  We 
call $\mathcal A$ the {\it ambient algebra} of $\Phi$.
\end{definition}

\section{The relatives of a Leonard system}

\medskip
\noindent A given Leonard system  can be modified in  several
ways to get a new Leonard system. For instance, 
let $\Phi$ 
 denote the Leonard system from 
(\ref{eq:ourstartingpt}),
 and let $\alpha$, $\alpha^*$, $\beta$,
$\beta^*$ denote scalars in $\K$ such that $\alpha \not=0$, $\alpha^*\not=0$.
Then
\begin{eqnarray}
(\alpha A+\beta I,
\alpha^* A^*+\beta^* I; E_i, E^*_i,i=0..d )
\label{eq:extendrel}
\end{eqnarray}
is a Leonard system in $\alg$.
Also, each of the following three sequences is a Leonard system
in $\mathcal A$.
%\begin{eqnarray}
\beast
\;\Phi^*&:=& (A^*,A; E^*_i,E_i,i=0..d),
%\label{eq:lsdualS99}
\\
\Phi^{\downarrow}&:=& (A,A^*; E_i,E^*_{d-i}, i=0..d),
%\label{eq:lsinvertS99}
\\
\Phi^{\Downarrow} 
&:=& (A, A^*; E_{d-i},E^*_i,i=0..d).
%\label{eq:lsdualinvertS99}
%\end{eqnarray}
\eeast
We refer to $\Phi^*$
(resp.  
 $\Phi^\downarrow$)   
(resp.  
 $\Phi^\Downarrow$) 
 as the  
{\it dual} 
(resp. {\it first inversion})
(resp.  {\it second inversion}) of  $\Phi$.
Viewing $*, \downarrow, \Downarrow$
as permutations on the set of all Leonard systems,
\begin{eqnarray}
&&\qquad \qquad \qquad  *^2 \;=\;  
\downarrow^2\;= \;
\Downarrow^2 \;=\;1,
\qquad \quad 
\label{eq:deightrelationsAS99}
\\
&&\Downarrow *\; 
=\;
* \downarrow,\qquad \qquad   
\downarrow *\; 
=\;
* \Downarrow,\qquad \qquad   
\downarrow \Downarrow \; = \;
\Downarrow \downarrow.
\qquad \quad 
\label{eq:deightrelationsBS99}
\end{eqnarray}
The group generated by symbols 
$*, \downarrow, \Downarrow $ subject to the relations
(\ref{eq:deightrelationsAS99}),
(\ref{eq:deightrelationsBS99})
is the dihedral group $D_4$.  
We recall $D_4$ is the group of symmetries of a square,
and has 8 elements.
Apparently $*, \downarrow, \Downarrow $ induce an action of 
 $D_4$ on the set of all Leonard systems.
Two Leonard systems will be called {\it relatives} whenever they
are in the same orbit of this $D_4$ action.
%Assuming $d\geq 1$ to avoid
%trivialities, the relatives of $\Phi$ are as follows:
The relatives of $\Phi$ are as follows:
\medskip

\centerline{
\begin{tabular}[t]{c|c}
        name &relative \\ \hline 
        $\Phi$ & $(A,A^*;E_i,E^*_i,i=0..d)$   \\ 
        $\Phi^\downarrow$ &
         $(A,A^*;E_i,E^*_{d-i},i=0..d)$   \\ 
        $\Phi^\Downarrow$ &
         $(A,A^*;E_{d-i},E^*_i,i=0..d)$   \\ 
        $\Phi^{\downarrow \Downarrow}$ &
         $(A,A^*;E_{d-i},E^*_{d-i},i=0..d)$   \\ 
	$\Phi^*$ &
        $(A^*,A;E^*_i,E_i,i=0..d)$   \\ 
        $\Phi^{\downarrow *}$ &
	 $(A^*,A;E^*_{d-i},E_i,i=0..d)$  
	 \\
        $\Phi^{\Downarrow *}$ &
	 $(A^*,A;E^*_i,    
       E_{d-i},i=0..d)$ \\ 
	$\Phi^{\downarrow \Downarrow *}$ &
	 $(A^*,A;E^*_{d-i}, E_{d-i},i=0..d)$    
	\end{tabular}}

\section{Leonard pairs and Leonard systems}

\medskip
In view of our comments in the previous section, when we discuss
a Leonard system we are often not interested in the orderings
of the primitive idempotents, we just care how the elements 
$A,A^*$ interact. This brings us back to the notion of a Leonard
pair.

\begin{definition} 
%\cite[Lem. 1.7]{LS99}
\label{def:lpabs}
\rm
Let $d$  denote  a  nonnegative
integer
 and let $\alg$ 
denote a $\K$-algebra isomorphic to 
$\hbox{Mat}_{d+1}(\K)$. 
By a {\it Leonard pair in 
 $\alg$} we mean an ordered pair $A,A^*$ which satisfies 
(i)--(iii) below.
\begin{enumerate}
\item Each of $A,A^*$ is a multiplicity-free element of $\mathcal A$.
\item There exists an idempotent sequence for $A,A^*$.
\item There exists a dual idempotent sequence for $A,A^*$.
\end{enumerate}
\end{definition}
By \cite[Lemma 1.7]{LS99}
the preceding definition of a Leonard pair is equivalent to 
the definition given in the Introduction. 

\medskip
\noindent Let $\Phi$ denote the Leonard system from 
(\ref{eq:ourstartingpt}). Then the pair $A,A^*$ from
that line forms a Leonard pair in $\mathcal A$.
We say this pair is  {\it associated} with $\Phi$.

\medskip
\noindent 
Each Leonard system is associated with a unique Leonard
pair. 
Let $A,A^*$ denote a Leonard pair.
By the {\it associate class} for $A,A^*$ we mean
the set of Leonard systems which are associated with
$A,A^*$.
By Definition
\ref{def:lpabs} this associate class contains
at least one Leonard system $\Phi$.
Apparently this associate class
contains $\Phi$,
$\Phi^\downarrow$, 
 $\Phi^\Downarrow$, 
$\Phi^{\downarrow \Downarrow}$
and no other Leonard systems.

\medskip
\noindent 
Let $A,A^*$ denote the Leonard pair from  
Definition
\ref{def:lpabs}.
Then the pair $A^*, A$ is a Leonard pair in $\mathcal A$.
We call this pair the {\it dual} of $A,A^*$.
We observe two Leonard systems are relatives if and only
if their associated Leonard pairs are equal or dual.

\section{Isomorphisms of Leonard pairs and Leonard systems}

\noindent We recall the notion of {\it isomorphism}
for Leonard pairs and Leonard systems.
We begin with a comment.

\begin{lemma}
\label{lem:uniquemap}
\cite[Corollary 3.2]{LS24}
Let $A,A^*$ denote the Leonard pair 
from  Definition \ref{def:lpabs}.
 Then $A$ and $A^*$ together generate
$\mathcal A$.
\end{lemma}

\medskip
\noindent 
Let $\Phi$ 
denote the Leonard system  from 
(\ref{eq:ourstartingpt})  
and let 
$\sigma :\alg \rightarrow {\mathcal A}'$ denote an isomorphism of
$\K$-algebras. We write 
$\Phi^{\sigma}:= 
(A^{\sigma},
A^{*\sigma}; 
 E_i^{\sigma},
E_i^{*\sigma}, i=0..d) $
and observe 
$\Phi^{\sigma}$ 
is a Leonard  system in ${\mathcal A }'$.

\begin{definition}
% \cite{LS99}
% \cite{TD00}, \cite{LS99}.
\label{def:isolsS99o}
\rm 
Let $\Phi$ 
and  
 $\Phi'$ 
denote Leonard systems over $\K$.
 By an {\it isomorphism of Leonard  systems}
 from $\Phi $ to $\Phi'$ we mean an isomorphism $\sigma$
 of $\K $-algebras
 from the ambient algebra of $\Phi$ to the ambient algebra of
$\Phi'$  such that
 $\Phi^\sigma = \Phi'$. 
By Lemma \ref{lem:uniquemap} there exists at most
one isomorphism of Leonard systems from 
 $\Phi $ to $\Phi'$. 
We say $\Phi $ and $\Phi'$
are {\it isomorphic} whenever this isomorphism
exists.
\end{definition}

\noindent 
We now consider the notion of isomorphism for Leonard pairs.

\medskip
\noindent 
Let $A,A^*$ denote the Leonard pair from Definition
\ref{def:lpabs} 
and let $\sigma: {\mathcal A}\rightarrow 
{\mathcal A}'$ denote an isomorphism of $\K$-algebras.
We observe the pair $A^\sigma, A^{*\sigma}$ is a Leonard
pair in ${\mathcal A}'$.

\begin{definition}
\rm 
Let $A,A^*$  
and  
 $B,B^*$ 
denote Leonard pairs over $\K$. 
 By an {\it isomorphism of Leonard pairs} 
 from $A,A^* $ to $B,B^*$ we mean an isomorphism $\sigma $
 of $\K $-algebras
from the ambient algebra of $A,A^*$ to the ambient algebra of
$B,B^*$ such that
 $A^\sigma =B$  and $A^{*\sigma}=B^*$.
By Lemma 
\ref{lem:uniquemap} there exists at most one 
isomorphism of Leonard pairs from $A,A^*$ to $B,B^*$. 
We say $A, A^*$ and $B,B^*$
are {\it isomorphic} whenever this isomorphism
exists.
\end{definition}

\noindent We have a comment.

\begin{lemma}
\label{lem:nonisoprov}
Let $A,A^*$ denote a Leonard pair and let $d$ denote
the diameter.
If $d\geq 1$ then 
the corresponding associate class contains
 four
Leonard systems and these are mutually nonisomorphic.
If 
$d= 0$ then the corresponding associate class contains a single
Leonard system.
\end{lemma}
\noindent {\it Proof:}
Let $\Phi$ denote a Leonard system associated
with $A,A^*$. Then the associate class of $\Phi$ contains
$\Phi$,
$\Phi^\downarrow$, 
 $\Phi^\Downarrow$, 
$\Phi^{\downarrow \Downarrow}$
and no other Leonard systems.
Suppose $d \geq 1$. Then    
$\Phi$,
$\Phi^\downarrow$, 
 $\Phi^\Downarrow$, 
$\Phi^{\downarrow \Downarrow}$
are mutually nonisomorphic; if not
the isomorphism  involved 
would stabilize each of $A,A^*$  
and is therefore the identity map
by Lemma \ref{lem:uniquemap}.
Suppose $d=0$. Then
$\Phi$,
$\Phi^\downarrow$, 
 $\Phi^\Downarrow$, 
$\Phi^{\downarrow \Downarrow}$
are identical by the construction.
 \hfill $\Box $\\

\noindent We finish this section with a remark.
Let $\mathcal A$ denote a matrix algebra over $\K$.
Let  $\sigma :\mathcal A \rightarrow 
\mathcal A$ denote  any map.
By the Skolem-Noether theorem 
\cite[Corollary 9.122]{JR},  
$\sigma $ is an isomorphism of $\K$-algebras
if and only if there exists an invertible $S \in \mathcal A$ such that
$X^\sigma = S X S^{-1}$ for all  $X \in  \mathcal A$.

\section{The adjacency relations}

\begin{definition}
\label{def:adj}
\rm 
Let $A,A^*$ denote the Leonard pair from Definition
\ref{def:lpabs}.
Consider the set consisting of the primitive idempotents of $A$.
We define a symmetric binary relation $\sim$
on this set.
Let  
$E_0,E_1,\ldots,
E_d$
denote an idempotent sequence for $A,A^*$.
For $0 \leq i, j \leq d$ we define $E_i\sim E_j$
whenever  $|i-j|=1$.
We call $\sim $  the
{\it first adjacency relation}
for $A,A^*$. We let $\approx $ denote the first adjacency
relation for the Leonard pair $A^*,A$ and call $\approx $
the 
{\it second adjacency relation}
for $A,A^*$.
\end{definition}

\noindent We make several observations. 

\begin{lemma}
\label{lem:evseqmeaning}
Let $A,A^*$ denote the Leonard pair from  
Definition
\ref{def:lpabs}.
Let $E_0, E_1, \ldots, E_d$
(resp. 
$E^*_0, E^*_1, $ $ \ldots, E^*_d$)
denote an ordering of the primitive idempotents
of $A$ 
(resp. $A^*$.)
 Then $E_0, E_1, \ldots, $ $ E_d$ is an idempotent
sequence for $A,A^*$ if and only  if
$E_0 \sim E_1 \sim \cdots \sim E_d$.
Moreover 
$E^*_0, E^*_1, $ $ \ldots, E^*_d$ 
is a dual idempotent sequence for $A,A^*$ if and only if  
$E^*_0 \approx E^*_1 \approx \cdots \approx E^*_d$.
\end{lemma}

\begin{lemma}
\label{lem:adjcheck}
Let $A,A^*$ denote the Leonard pair 
from
Definition
\ref{def:lpabs}. Let $E$ and $F$ denote primitive idempotents
of $A$. Then the following are equivalent: (i) $E\sim F$;
(ii)
  $E\not=F$ and $EA^*F\not=0$;
 (iii) $E\not=F$ and $FA^*E\not=0$.
Let $E^*$ and $F^*$ denote primitive idempotents
of $A^*$. Then 
 the following are equivalent: (i) $E^* \approx F^*$;
(ii)
  $E^*\not=F^*$ and $E^*AF^*\not=0$;
 (iii) $E^*\not=F^*$ and $F^*AE^*\not=0$.
\end{lemma}

\section{The eigenvalue sequences}

\begin{definition}
\label{def:evseq}
\rm
Let $\Phi$ denote the Leonard system from  
(\ref{eq:ourstartingpt}).
For $0 \leq i \leq d$  
we let $\theta_i $ (resp. $\theta^*_i$) denote the eigenvalue
of $A$ (resp. $A^*$) associated with $E_i$ (resp. $E^*_i$.)
We call  $\theta_0, \theta_1, \ldots, \theta_d$  the 
{\it eigenvalue sequence} of $\Phi$.
We call  $\theta^*_0, \theta^*_1, \ldots, \theta^*_d$  the 
{\it dual eigenvalue sequence} of $\Phi$. We observe 
 $\theta_0, \theta_1, \ldots, \theta_d$ are mutually distinct
 and contained in $\K$. Similarly
  $\theta^*_0, \theta^*_1, \ldots, \theta^*_d$  
 are mutually distinct
 and contained in $\K$. 
\end{definition}

\begin{definition}
\label{def:evseqlp}
\rm 
Let $A,A^*$ denote a Leonard pair.
By an {\it eigenvalue sequence} for this pair, we mean
the eigenvalue sequence for an associated Leonard system.
We remark that if
$\theta_0, \theta_1, \ldots,\theta_d$ is an eigenvalue
sequence for $A,A^*$ then
so is  $\theta_d, \theta_{d-1}, \ldots,\theta_0$ and
$A,A^*$ has no further eigenvalue sequence.
By a {\it  dual eigenvalue sequence} for $A,A^*$
 we mean an eigenvalue sequence for the  
Leonard pair $A^*,A$.
\end{definition}

\noindent We will use the following results.

\begin{lemma}
\label{lem:vague}
Let $d$ denote a nonnegative integer and let $A,A^*$ denote
a Leonard pair in 
$\hbox{Mat}_{d+1}(\K)$.
%We make the following assumptions about the entries of $A,A^*$.
Assume (i) $A$ is lower triangular; and (ii) 
$A^*_{ij}=0$
if $j-i>1$, $(0\leq i,j\leq d)$.
Then the sequence of 
diagonal entries $A_{00}, A_{11},\ldots, A_{dd}$ is an eigenvalue
sequence for $A,A^*$. Moreover $A^*_{j-1,j}\not=0$ 
for $1\leq j \leq d$.
\end{lemma}
\noindent {\it Proof:}
We assume the pair $A,A^*$ is a Leonard pair so
$A$ is multiplicity-free.
We assume $A$ is lower triangular so the sequence of
diagonal entries 
$A_{00}, A_{11},\ldots, A_{dd}$ gives an ordering of
the eigenvalues of $A$.
We show this sequence is an eigenvalue sequence for $A,A^*$.
For $0 \leq i \leq d$ 
let $E_i$ denote the primitive idempotent of $A$ associated
with the eigenvalue $A_{ii}$.
We show $E_{j-1}\sim E_j$ for $1 \leq j\leq d$.
This will follow once we 
show 
\begin{eqnarray}
E_i \not\sim E_j 
\qquad \hbox{if} \quad  j-i > 1 \qquad \qquad 
(0 \leq i,j\leq d).
\label{eq:part1}
\end{eqnarray}
We abbreviate $V=
\K^{d+1}$.
For $0 \leq r \leq d$ let $V_r$ denote the subspace of $V$
consisting of those vectors which have 0 in coordinates
$0,1,\ldots, r-1$. The matrix $A$ is lower triangular so
$AV_r \subseteq V_r$. The restriction of $A$ to
$V_r$ has eigenvalues $A_{rr}, \ldots, A_{dd}$ so
 $V_r=E_rV+\cdots + E_dV$. Apparently $E_rV \subseteq V_r$
 and moreover each of $E_0, \ldots, E_{r-1}$ vanishes on $V_r$.
From our assumption about  $A^*$ we find $A^*V_r \subseteq V_{r-1}$
for $1 \leq r \leq d$. 
Let $i,j$ denote integers $(0 \leq i,j\leq d)$ and assume
$j-i>1$. From our above comments we find
\beast
E_iA^*E_jV \subseteq E_iA^*V_j \subseteq E_iV_{j-1} = 0.
\eeast
Apparently  $E_iA^*E_jV=0$ so 
$E_iA^*E_j=0$. Now  $E_i {\not \sim} E_j $ by Lemma
\ref{lem:adjcheck}.
We now have 
(\ref{eq:part1}) and it follows
$E_{j-1}\sim E_j$ for $1 \leq j \leq d$.
Applying Lemma
\ref{lem:evseqmeaning} we find
$E_0, E_1, \ldots, E_d$ is an idempotent sequence for
$A,A^*$. Now
$A_{00}, A_{11},\ldots, A_{dd}$ is an eigenvalue
sequence for $A,A^*$ by Definition
\ref{def:evseqlp}.
To finish the proof we show $A^*_{j-1,j}\not=0$ for
$1 \leq j\leq d$.
Let $j$ be given and suppose $A^*_{j-1,j}=0$. Then
$A^*V_j\subseteq V_j$.
We mentioned earlier that $AV_j \subseteq V_j$. 
The matrices 
$A$ and $A^*$  together 
generate
$\hbox{Mat}_{d+1}(\K)$
by Lemma
\ref{lem:uniquemap}
so $XV_j \subseteq V_j$ for
all 
$X \in \hbox{Mat}_{d+1}(\K)$. The space 
$V$ is irreducible as a module for
$\hbox{Mat}_{d+1}(\K)$,  so 
 $V_j =0$ or $V_j=V$.
From the definition of $V_j$ and since $1 \leq j \leq d$ we
find $V_j \not=0$ and $V_j\not=V$. 
This is a contradiction and we conclude 
 $A^*_{j-1,j}\not=0$.
\hfill $\Box $\\

\begin{lemma}
\label{lem:vague2}
Let $d$ denote a nonnegative integer and let $A,A^*$ denote
a Leonard pair in
$\hbox{Mat}_{d+1}(\K)$.
Assume (i) $A$ is upper triangular; and (ii) $A^*_{ij}=0$
if $i-j>1$, $(0\leq i,j\leq d)$.
Then the sequence of 
diagonal entries $A_{00}, A_{11},\ldots, A_{dd}$ is an eigenvalue
sequence for $A,A^*$. Moreover $A^*_{i,i-1}\not=0$ for 
$1\leq i \leq d$.
\end{lemma}
\noindent {\it Proof:}
Using Definition
\ref{def:lpabs}
we find
$A^t, A^{*t}$ is a Leonard pair  in
$\hbox{Mat}_{d+1}(\K)$,
where $t$ denotes transpose.
To obtain the result apply Lemma
\ref{lem:vague} to this pair.
\hfill $\Box $\\

\noindent We give a corollary to Lemma
\ref{lem:vague} and
Lemma \ref{lem:vague2}.
 In order to state it we make a definition.

\begin{definition}
\label{def:lbubmat}
\rm 
Let $d$ denote a nonnegative integer and let
$A$ denote a matrix in 
$\hbox{Mat}_{d+1}(\K)$. 
We say $A$ is {\it lower bidiagonal} whenever
each nonzero entry  lies  on either the diagonal or
the subdiagonal. We say $A$ is {\it upper bidiagonal}
whenever the transpose of $A$ is lower bidiagonal.
\end{definition}

\begin{corollary}
\label{lem:basicdata}
Let $d$ denote a nonnegative integer and let 
$A,A^*$ denote a Leonard pair in 
$\hbox{Mat}_{d+1}(\K)$. Assume
$A$ is lower bidiagonal and $A^*$ is upper bidiagonal.
Then (i)--(iv) hold below. 
\begin{enumerate}
\item The sequence $A_{00}, A_{11}, \ldots, A_{dd}$
is an eigenvalue sequence for $A,A^*$.
\item $A_{i,i-1}\not=0$ for $1 \leq i\leq d$.
\item The sequence $A^*_{00}, A^*_{11}, \ldots, A^*_{dd}$
is a dual eigenvalue sequence for $A,A^*$.
\item $A^*_{i-1,i}\not=0$ for $1 \leq i\leq d$.
\end{enumerate}
\end{corollary}
\noindent {\it Proof:}
(i),(iv) Apply Lemma 
\ref{lem:vague} to $A,A^*$.
\\
\noindent (ii),(iii) Apply Lemma 
\ref{lem:vague2} to the Leonard pair $A^*,A$.
\hfill $\Box $\\

\noindent The following fact
may seem 
intuitively clear from
Definition
\ref{def:lpabs}, but strictly 
speaking it requires proof.

\begin{corollary}
\label{cor:tdeig}
Let $d$ denote a nonnegative integer and let
$A,A^*$ denote a Leonard pair in 
$\hbox{Mat}_{d+1}(\K)$. Assume
$A$ is tridiagonal and $A^*$ is diagonal.
Then (i), (ii) hold below.
\begin{enumerate}
\item $A$ is irreducible.
\item The sequence  
 $A^*_{00}, A^*_{11}, \ldots, A^*_{dd}$ is
 a dual eigenvalue sequence for $A,A^*$.
\end{enumerate}
\end{corollary}
\noindent {\it Proof:}
(i) Applying Lemma
\ref{lem:vague2}
to the Leonard pair $A^*,A$ we find $A_{i,i-1}\not=0$ for $1 \leq i \leq d$.
Applying Lemma
\ref{lem:vague}
to $A^*,A$ we find $A_{i-1,i}\not=0$ for $1 \leq i \leq d$.
\\
\noindent (ii) Apply Lemma
\ref{lem:vague} to the Leonard pair $A^*,A$.
\hfill $\Box $\\

\section{The split sequences}

\noindent 
In 
Definition \ref{def:evseq}
we defined the eigenvalue sequence
and the dual eigenvalue sequence of a Leonard system.
There are two more parameter sequences of interest to us.
In order to define these, we review some results from 
  \cite{TD00},
\cite{LS24},
\cite{LS99}. 
Let $\Phi$ denote the Leonard system
in 
(\ref{eq:ourstartingpt})
and let $V$ denote an irreducible
 $\mathcal A$-module.
For $0 \leq i \leq d$ we define
\begin{equation}
U_i = 
(E^*_0V + E^*_1V + \cdots + E^*_iV)\cap (E_iV + E_{i+1}V + \cdots + E_dV).
\label{eq:defui}
\end{equation}
We showed in \cite[Lemma 3.8]{LS99}
that each of $U_0, U_1, \ldots, U_d$ has dimension 1, and that
\begin{equation}
V = U_0 + U_1 + \cdots + U_d \qquad \qquad (\hbox{direct sum}).
\label{eq:splitdec}
\end{equation}
The elements $A$ and $A^*$ act on the $U_i$ as follows.
By \cite[Lemma 3.9]{LS99}, both
\begin{eqnarray}
(A-\theta_i I)U_i &=& U_{i+1} \quad (0 \leq i \leq d-1),
\qquad (A-\theta_d I)U_d = 0,
\label{eq:raise}
\\
(A^*-\theta^*_i I)U_i &=& U_{i-1} \quad (1 \leq i \leq d),
\qquad (A^*-\theta^*_0 I)U_0 = 0,
\label{eq:lower}
\end{eqnarray}
where the $\theta_i, \theta^*_i$ are from
Definition \ref{def:evseq}. 
Pick an integer $i$ $(1 \leq i \leq d)$. By
(\ref{eq:lower}) we find 
$(A^*-\theta^*_i I)U_i = U_{i-1}$ and by
(\ref{eq:raise}) we find 
$(A-\theta_{i-1} I)U_{i-1} = U_i$. Apparently $U_i$
is an eigenspace for
$(A-\theta_{i-1}I)(A^*-\theta^*_i I)$
and the corresponding eigenvalue is a nonzero element of $\K$.
We denote this 
eigenvalue by $\varphi_i$.
We display a basis for $V$ which illuminates the significance
of  $\varphi_i$.  
Setting $i=0$ in 
(\ref{eq:defui}) 
we 
find $U_0=E^*_0V$. Combining this with
(\ref{eq:raise}) we find  
\begin{equation}
U_i = (A-\theta_{i-1} I)\cdots (A - \theta_1 I)(A-\theta_0 I )E^*_0V
\qquad \qquad (0 \leq i\leq d).
\label{eq:uialt}
\end{equation}
Let $\eta^*_0$ denote a nonzero vector in $E^*_0V$. 
From 
(\ref{eq:uialt}) we find that for 
$0 \leq i \leq d$ the vector
$(A-\theta_{i-1}I)\cdots (A-\theta_{0}I)\eta^*_0 $ is a basis for
$U_i$. From this and 
(\ref{eq:splitdec}) we find the sequence
\begin{equation}
(A-\theta_{i-1} I)\cdots (A-\theta_1 I)(A-\theta_{0}I)\eta^*_0  \qquad \qquad
(0 \leq i \leq d)
\label{eq:basis1}
\label{eq:sp}
\end{equation}
is a basis for $V$. 
With respect to this basis
the matrices representing $A$ and $A^*$ are 
\begin{equation}
\left(
\begin{array}{c c c c c c}
\theta_0 & & & & & {\bf 0} \\
1 & \theta_1 &  & & & \\
& 1 & \theta_2 &  & & \\
& & \cdot & \cdot &  &  \\
& & & \cdot & \cdot &  \\
{\bf 0}& & & & 1 & \theta_d
\end{array}
\right),
\qquad  \quad 
%A^{*\heartsuit} = 
%A^{*\flat} = 
\left(
\begin{array}{c c c c c c}
\theta^*_0 &\varphi_1 & & & & {\bf 0} \\
 & \theta^*_1 & \varphi_2 & & & \\
&  & \theta^*_2 & \cdot & & \\
& &  & \cdot & \cdot &  \\
& & &  & \cdot & \varphi_d \\
{\bf 0}& & & &  & \theta^*_d
\end{array}
\right)
\label{eq:matrepaastar}
\end{equation}
respectively. 
 We call 
 the sequence $\varphi_1, \varphi_2,
\ldots, \varphi_d$  the {\it first split sequence} of 
$\Phi$. We let $\phi_1, \phi_2, \ldots, \phi_d$ denote
the first split sequence for $\Phi^{\Downarrow}$ and
call  this the {\it  second split sequence} of $\Phi$. For notational
convenience we define $\varphi_0=0$, $\varphi_{d+1} = 0$,
 $\phi_0=0$, $\phi_{d+1} = 0$.

\section{A classification of Leonard systems}

We recall our classification of Leonard systems.  

\begin{theorem} \cite[Theorem 1.9]{LS99}
\label{thm:classls}
Let 
$d$ denote a nonnegative integer  
and let 
%\begin{eqnarray}
\beast
&&\theta_0, \theta_1, \ldots, \theta_d; \qquad \qquad \; 
\theta^*_0, \theta^*_1, \ldots, \theta^*_d; 
%\label{eq:paramlist1S99}
\\
&&\varphi_1, \varphi_2, \ldots, \varphi_d;  \qquad \qquad 
\phi_1, \phi_2, \ldots, \phi_d \qquad \quad
%\label{eq:paramlist2S99}
%\end{eqnarray}
\eeast
denote scalars in $\K$. 
Then there exists  a Leonard system $\Phi$ over $\K$  with 
eigenvalue sequence $\theta_0, \theta_1, \ldots, \theta_d$, 
dual eigenvalue sequence  
$\theta^*_0, \theta^*_1, \ldots, \theta^*_d$, first split sequence
$ \varphi_1, \varphi_2, \ldots, \varphi_d $ and second split sequence
$\phi_1, \phi_2, \ldots, \phi_d$ if and only if 
(i)--(v) hold below.
\begin{enumerate}
\item $ \varphi_i \not=0, \qquad \phi_i\not=0 \qquad \qquad \qquad\qquad (1 \leq i \leq d)$.
\item $ \theta_i\not=\theta_j,\qquad  \theta^*_i\not=\theta^*_j\qquad $
if $\;\;i\not=j,\qquad \qquad \qquad (0 \leq i,j\leq d)$.
\item $ {\displaystyle{ \varphi_i = \phi_1 \sum_{h=0}^{i-1}
{{\theta_h-\theta_{d-h}}\over{\theta_0-\theta_d}} 
\;+\;(\theta^*_i-\theta^*_0)(\theta_{i-1}-\theta_d) \qquad \;\;(1 \leq i \leq d)}}$.
\item $ {\displaystyle{ \phi_i = \varphi_1 \sum_{h=0}^{i-1}
{{\theta_h-\theta_{d-h}}\over{\theta_0-\theta_d}} 
\;+\;(\theta^*_i-\theta^*_0)(\theta_{d-i+1}-\theta_0) \qquad (1 \leq i \leq d)}}$.
\item The expressions
\beast
%\begin{equation}
{{\theta_{i-2}-\theta_{i+1}}\over {\theta_{i-1}-\theta_i}},\qquad \qquad  
 {{\theta^*_{i-2}-\theta^*_{i+1}}\over {\theta^*_{i-1}-\theta^*_i}} 
 \qquad  \qquad 
%\label{eq:betaplusone}
%\end{equation} 
\eeast
are equal and independent of $i$ for $\;2\leq i \leq d-1$.  
\end{enumerate}
Moreover, if (i)--(v) hold 
above then $\Phi$ is unique up to isomorphism of Leonard systems.
\end{theorem}
\noindent We view Theorem
\ref{thm:classls} 
as a linear algebraic version
of a theorem of D. Leonard
 \cite{Leon}, \cite[p. 260]{BanIto}. This is discussed in \cite{LS99}.

\section{The notion of a parameter array}

\noindent In view of Theorem
\ref{thm:classls} we make the following definition. 

\begin{definition}
\rm 
Let $d$ denote a nonnegative integer.
By a {\it parameter array over $\K$ with diameter $d$},
we mean a sequence  
$(\theta_i, \theta^*_i, i=0..d;  \varphi_j, \phi_j, j=1..d)$
of scalars taken from $\K$ which satisfy conditions (i)--(v) 
in Theorem
\ref{thm:classls}. 
\end{definition}

\noindent 
We give several examples of a parameter array.

\begin{example}
\label{ex:pa1}
\rm
Let $d$ denote a nonnegative integer and 
consider the following scalars
in $\K$.
%\begin{eqnarray}
\beast
\theta_i&=&d-2i, \qquad \qquad \theta^*_i=d-2i \qquad \qquad
(0 \leq i \leq d),
%\label{eq:ththf}
\\
\varphi_i&=&-2i(d-i+1), \qquad \qquad \phi_i=2i(d-i+1)
\qquad \qquad (1 \leq i \leq d).
%\label{eq:split2f}
\eeast
%\end{eqnarray}
To avoid degenerate situations,
we assume the characterisic of $\K$ is zero or an odd prime
greater than $d$.
Then the sequence 
$(\theta_i, \theta^*_i, i=0..d;  \varphi_j, \phi_j, j=1..d)$
is a parameter array over $\K$.
\end{example}
\noindent {\it Proof:}
The sequence 
$(\theta_i, \theta^*_i, i=0..d;  \varphi_j, \phi_j, j=1..d)$
satisfies Theorem
\ref{thm:classls}(i)--(v) so 
this sequence
is a parameter array over $\K$.
\hfill $\Box $\\

\begin{example}
\label{ex:pa2}
\label{ex:qracah}
\rm
Let $d$ denote a nonnegative integer.
Let 
 $q, s, s^*, r_1, r_2$ denote nonzero scalars in 
$\K$ such that $r_1r_2=ss^*q^{d+1}$.
Assume
none of
$q^i,
r_1q^i, r_2q^i,
s^*q^i/r_1,$ $ s^*q^i/r_2$
is equal to $1$ for $1 \leq i \leq d$ and 
that neither of $sq^i, s^*q^i$ is equal to $1$ for $2 \leq i \leq 2d$.
Define
\beast
\theta_i &=& q^{-i}+ sq^{i+1}, 
\qquad \qquad 
\theta^*_i = 
 q^{-i}+ s^*q^{i+1}
\eeast
for $0 \leq i \leq d$ and
\beast
\varphi_i &=& q^{1-2i}(1-q^i)(1-q^{i-d-1})(1-r_1q^i)(1-r_2q^i),
\\
\phi_i &=& q^{1-2i}(1-q^i)(1-q^{i-d-1})(r_1-s^*q^i)(r_2-s^*q^i)/s^*
\eeast
for $1 \leq i \leq d$.
Then the sequence
$(\theta_i, \theta^*_i, i=0..d;  \varphi_j, \phi_j, j=1..d)$
is a parameter array over $\K$.
\end{example}
\noindent {\it Proof:}
The sequence
$(\theta_i, \theta^*_i, i=0..d;  \varphi_j, \phi_j, j=1..d)$
satisfies  Theorem 
\ref{thm:classls}(i)--(v)  
so this sequence
is a parameter array over $\K$.
\hfill $\Box $\\

\section{Parameter arrays and Leonard systems}

In this section  we discuss the relationship between 
parameter arrays and
Leonard systems.

\begin{definition}
\label{def:pals}
\rm
Let $\Phi$ denote a Leonard system over
$\K$, with 
eigenvalue sequence $\theta_0, \theta_1, \ldots$, $\theta_d$, 
dual eigenvalue sequence  
$\theta^*_0, \theta^*_1, \ldots, \theta^*_d$, first split sequence
$ \varphi_1, \varphi_2, \ldots, \varphi_d $, and second split sequence
$\phi_1, \phi_2, \ldots, \phi_d$.
By 
 Theorem
\ref{thm:classls}
the sequence
$(\theta_i, \theta^*_i, i=0..d;  \varphi_j, \phi_j, j=1..d)$
is a parameter array over $\K$.
We call this array the {\it parameter array of} $\Phi$.   
\end{definition}

\noindent
We remark that by
 Theorem
\ref{thm:classls}
the map which sends a given Leonard system 
 to its parameter array induces a bijection
from the set of isomorphism
classes of
Leonard systems 
over $\K$ to the set of parameter
arrays over $\K$.

\medskip
\noindent  
Earlier we discussed several ways to modify a given Leonard system
to get a new Leonard system. We now consider how these modifications
affect 
the corresponding  
parameter array.

\begin{lemma}
Let $\Phi$ denote the Leonard system from  
(\ref{eq:ourstartingpt})
 and let  
$(\theta_i, \theta^*_i, i=0..d;  \varphi_j, \phi_j, j=1..d)$
denote the corresponding parameter array.
 Let $\alpha$, $\alpha^*$, $\beta$,
$\beta^*$ denote scalars in $\K$ such that $\alpha \not=0$, $\alpha^*\not=0$.
Then the Leonard system
(\ref{eq:extendrel}) has parameter array
\beast
(\alpha \theta_i+\beta, \alpha^*\theta^*_i+\beta^*, i=0..d; 
\alpha \alpha^*\varphi_j, \alpha \alpha^* \phi_j, j=1..d).
\eeast
\end{lemma}
\noindent {\it Proof:}
Routine.
\hfill $\Box $\\

\begin{lemma}
\label{thm:phimod}
\cite[Theorem 1.11]{LS99}
Let $\Phi$ denote a Leonard system and let  
$p=(\theta_i, \theta^*_i, i=0..d;  \varphi_j, \phi_j,$ $ j=1..d)$
denote the corresponding parameter array. Then (i)--(iii) hold below.
\begin{enumerate}
\item The parameter array of $\Phi^*$ is $p^*$ where
$p^*:=(\theta^*_i, \theta_i, i=0..d;  \varphi_j, \phi_{d-j+1}, j=1..d).$
\item 
The parameter array of $\Phi^\downarrow$ is $p^\downarrow$
where
$p^\downarrow :=(\theta_i, \theta^*_{d-i}, i=0..d;  \phi_{d-j+1}, \varphi_{d-j+1}, j=1..d).
$
\item 
The parameter array of $\Phi^\Downarrow$ is
$p^\Downarrow $ where
$p^\Downarrow:=
(\theta_{d-i}, \theta^*_i, i=0..d;  \phi_{j}, \varphi_{j}, j=1..d).
$
\end{enumerate}
\end{lemma}

\noindent The following equations will be useful.

\begin{corollary}
\label{cor:d4pa}
Let $d$ denote a positive integer and let
$(\theta_i, \theta^*_i, i=0..d;  \varphi_j, \phi_j, j=1..d)$
denote a parameter array over $\K$. Then (i)--(iii) hold below.
\begin{enumerate}
\item 
$ {\displaystyle{
{{\theta_i-\theta_{d-i}}\over{\theta_0-\theta_d}} 
={{\theta^*_i-\theta^*_{d-i}}\over{\theta^*_0-\theta^*_d}} 
\qquad \qquad \qquad (0 \leq i\leq d)
}}$.
\item $ {\displaystyle{ \varphi_i = \phi_d \sum_{h=0}^{i-1}
{{\theta_h-\theta_{d-h}}\over{\theta_0-\theta_d}} 
\;+\;(\theta_i-\theta_0)(\theta^*_{i-1}-\theta^*_d) 
\qquad \;\;(1 \leq i \leq d)}}$.
\item $ {\displaystyle{ \phi_i = \varphi_d \sum_{h=0}^{i-1}
{{\theta_h-\theta_{d-h}}\over{\theta_0-\theta_d}} 
\;+\;(\theta_{d-i}-\theta_d)(\theta^*_{i-1}-\theta^*_d) \qquad (1 \leq i \leq d)}}$.
\end{enumerate}
\end{corollary}
\noindent {\it Proof:}
Each of (i)--(iii) is an algebraic consequence of
the  conditions in Theorem \ref{thm:classls}.
Below we give a more intuitive proof using  
Lemma \ref{thm:phimod}. Let $\Phi$ denote a Leonard system
over $\K$ which has the given parameter array.
\\
\noindent (i) Applying Theorem 
\ref{thm:classls}(iv) to $\Phi^*$ and using Lemma
\ref{thm:phimod}(i) we obtain 
\begin{eqnarray}
 \phi_{d-i+1} = \varphi_1 \sum_{h=0}^{i-1}
{{\theta^*_h-\theta^*_{d-h}}\over{\theta^*_0-\theta^*_d}} 
\;+\;(\theta_i-\theta_0)(\theta^*_{d-i+1}-\theta^*_0)
\label{eq:varphiadj}
\end{eqnarray}
for $1 \leq i \leq d$.
To finish the proof, in 
(\ref{eq:varphiadj})
replace 
$i$ by
$d-i+1$ and compare the result with 
Theorem
\ref{thm:classls}(iv).
\\
\noindent (ii) Apply Theorem
\ref{thm:classls}(iii) to $\Phi^*$ and simplify the result using
(i) above and Lemma
\ref{thm:phimod}(i).
\\
\noindent (iii) Apply (ii) above to $\Phi^\Downarrow$ and use
Lemma
\ref{thm:phimod}(iii).
\hfill $\Box $\\

\section{The parameter arrays of a Leonard pair}

\noindent In this section we define the notion
of a parameter array for a Leonard pair.

\medskip
\noindent
\begin{definition}
\label{def:palp}
\rm
Let $A,A^*$ denote a Leonard pair.
By a {\it parameter array} of $A,A^*$ we mean
the parameter array of an associated
Leonard system.
\end{definition}
\noindent 

\noindent The parameter arrays of a Leonard pair are related 
as follows.

\begin{lemma}
\label{lem:palprel}
Let $A,A^*$ denote the Leonard pair from
Definition
\ref{def:lpabs}.
Let 
$p=(\theta_i, \theta^*_i, i=0..d; $ $ \varphi_j, \phi_j, j=1..d)$
denote a parameter array of $A,A^*$. Then the following are parameter
arrays of $A,A^*$.
\beast
&&p=(\theta_i, \theta^*_i, i=0..d;  \varphi_j, \phi_j, j=1..d),
\\
&&p^\downarrow = (\theta_i, \theta^*_{d-i}, i=0..d;  \phi_{d-j+1}, \varphi_{d-j+1}, j=1..d),
\\
&&p^\Downarrow=(\theta_{d-i}, \theta^*_i, i=0..d;  \phi_j, \varphi_j,
j=1..d),
\\
&&p^{\downarrow \Downarrow}=(\theta_{d-i}, \theta^*_{d-i}, i=0..d;  \varphi_{d-j+1}, \phi_{d-j+1}, j=1..d).
\eeast
The Leonard pair $A,A^*$ has no further parameter arrays.
\end{lemma}
\noindent {\it Proof:}
By Definition
\ref{def:palp} there exists
a Leonard system $\Phi$  which is
associated with $A,A^*$ and which has parameter array $p$.
The above  sequences are the parameter arrays  for 
$\Phi$, $\Phi^\downarrow$, $\Phi^\Downarrow$,
$\Phi^{\downarrow \Downarrow}$ and these are the Leonard systems
associated with $A,A^*$.
\hfill $\Box $\\

\begin{corollary}
\label{cor:numberpa}
Let 
$A,A^*$ 
denote the Leonard pair from 
Definition
\ref{def:lpabs}.
Then $A,A^*$ has exactly four parameter arrays if
$d\geq 1$ and a unique parameter array if $d=0$.
\end{corollary}
\noindent {\it Proof:}
Referring to Lemma
\ref{lem:palprel},
 the parameter arrays
$p, p^\downarrow, p^\Downarrow, p^{\downarrow \Downarrow}$
are mutually distinct if $d\geq 1$ 
and identical if $d=0$.
\hfill $\Box $\\

\noindent 
We have a comment.

\begin{lemma}
\label{lem:parrayforlp}
Let $A,A^*$ denote a Leonard pair over $\K$ and let
$B,B^*$ denote a Leonard pair over $\K$.
These pairs are isomorphic if and only if 
they share a 
parameter array.
In this case the set of parameter arrays for $A,A^*$ coincides
with the set of parameter arrays for $B,B^*$.
\end{lemma}
\noindent {\it Proof:}
Suppose $A,A^*$ and $B,B^*$ share a parameter array
$p$. By Definition 
\ref{def:palp} there exists a Leonard system $\Phi$ which is
associated with $A,A^*$ and which has parameter array $p$.
Similarly 
there exists a Leonard system $\Phi'$ which is
associated with $B,B^*$ and which has parameter array $p$.
Observe $\Phi, \Phi'$ are isomorphic since they have 
the same parameter array.
Observe the isomorphism involved is an isomorphism of Leonard
pairs from $A,A^*$ to $B,B^*$. Apparently $A,A^*$ and $B,B^*$
are isomorphic. The remaining claims of the lemma are clear.
\hfill $\Box $\\

\section{The $LB$-$UB$ canonical form; preliminaries}

We now turn our attention to the $LB$-$UB$ canonical form.
We begin with some comments.

\begin{definition}
\label{lem:lbubbasis}
\rm
Let $\Phi$ denote the Leonard system from 
(\ref{eq:ourstartingpt}) and let 
$V$ denote an irreducible $\mathcal A$-module.
By a {\it $\Phi$-$LB$-$UB$ basis} for $V$ we mean
a sequence of the form
(\ref{eq:basis1}),
where $\theta_0, \theta_1, \ldots, \theta_d$ denotes the eigenvalue
sequence for $\Phi$ and 
$\eta^*_0$ denotes a nonzero vector in $E^*_0V$. 
\end{definition}

\begin{lemma}
\label{lem:lbubcheck}
Let $\Phi$ denote the Leonard system from 
(\ref{eq:ourstartingpt}). Let 
$\theta_0, \theta_1, \ldots, \theta_d$ denote the 
eigenvalue sequence for $\Phi$.
Let $V$ denote an 
irreducible $\mathcal A$-module and
let $v_0, v_1, \ldots, v_d$ denote a sequence of
vectors in $V$, not all zero. Then 
this sequence is a $\Phi$-$LB$-$UB$ basis for $V$
if and only if both (i)
$v_0 \in E^*_0V$; and 
(ii) $Av_i = \theta_iv_i + v_{i+1}$ for $0 \leq i \leq d-1$.
\end{lemma}
\noindent {\it Proof:}
Routine.
\hfill $\Box $\\

\begin{definition}
\label{def:natcon}
\rm
Let $\Phi$ denote the Leonard system from 
(\ref{eq:ourstartingpt}). We define
a map $\natural : {\mathcal A}\rightarrow 
\hbox{Mat}_{d+1}(\K)$ as follows.
Let 
$V$ denote an irreducible $\mathcal A$-module.
For all $X \in {\mathcal  A}$ we let $X^\natural $ denote
the matrix 
in 
$\hbox{Mat}_{d+1}(\K)$
which represents $X$ with respect to 
a $\Phi$-$LB$-$UB$ basis for $V$.
We observe $\natural : {\mathcal A} \rightarrow 
\hbox{Mat}_{d+1}(\K)$ is an isomorphism of
$\K$-algebras. We call $\natural $ the
{\it $LB$-$UB$ canonical map} for $\Phi$.
\end{definition}

\noindent Before proceeding we introduce some notation.

\begin{definition}
\label{def:LB}
\rm
Consider the set of all parameter arrays
over $\K$. 
We define two functions on this set. 
We call these functions 
$L$ and $U$.
Let 
$p=(\theta_i, \theta^*_i, i=0..d;  \varphi_j, \phi_j, j=1..d)$
denote a parameter array over $\K$.
The images $p^{L}$ and $p^U$ are the following matrices in
$\hbox{Mat}_{d+1}(\K)$.
%\begin{equation}
\beast
p^L=\left(
\begin{array}{c c c c c c}
\theta_0 & & & & & {\bf 0} \\
1 & \theta_1 &  & & & \\
& 1 & \theta_2 &  & & \\
& & \cdot & \cdot &  &  \\
& & & \cdot & \cdot &  \\
{\bf 0}& & & & 1 & \theta_d
\end{array}
\right),
\qquad  \quad 
p^U=
\left(
\begin{array}{c c c c c c}
\theta^*_0 &\varphi_1 & & & & {\bf 0} \\
 & \theta^*_1 & \varphi_2 & & & \\
&  & \theta^*_2 & \cdot & & \\
& &  & \cdot & \cdot &  \\
& & &  & \cdot & \varphi_d \\
{\bf 0}& & & &  & \theta^*_d
\end{array}
\right).
%\label{eq:plu}
%\end{equation}
\eeast
\end{definition}

\begin{lemma}
%\cite[Thm. 11.2]{LS24}
\label{lem:lbubobv}
Let $\Phi$ denote the 
Leonard system
from (\ref{eq:ourstartingpt}).
Let $\natural $ denote the $LB$-$UB$
canonical map for $\Phi$, from Definition
\ref{def:natcon}.
Then
$A^\natural=p^L$ and $A^{*\natural}=p^U$, where 
$p$ denotes the parameter array for $\Phi$.
\end{lemma}
\noindent {\it Proof:}
Write
$p=(\theta_i, \theta^*_i, i=0..d;  \varphi_j, \phi_j, j=1..d)$.
Each of $A^\natural, p^L$ is equal to the matrix
on the left in 
(\ref{eq:matrepaastar}) so
$A^\natural=p^L$.
Each of 
$A^{*\natural}, p^U$ is equal to the matrix
on the right in 
(\ref{eq:matrepaastar}) so
$A^{*\natural}=p^U$.
\hfill $\Box $\\

\section{The $LB$-$UB$ canonical form for Leonard systems} 

\noindent 
In this section we introduce the 
$LB$-$UB$ canonical form for Leonard systems.
We 
define what it means for a given Leonard system
 to be in $LB$-$UB$ canonical
form. We describe the 
  Leonard systems which are in 
$LB$-$UB$ canonical form.
 We  show 
every Leonard system is isomorphic to 
a unique Leonard system 
which is in $LB$-$UB$ canonical form.

\begin{definition}
\label{def:lbub}
\rm 
Let $\Phi$ denote the Leonard  system from 
(\ref{eq:ourstartingpt}). 
Let $\theta_0, \theta_1, \ldots, \theta_d$ (resp.  
 $\theta^*_0, \theta^*_1, \ldots, \theta^*_d$)
denote the eigenvalue sequence (resp. dual eigenvalue sequence) of
$\Phi$.
We say $\Phi$ is in 
{\it $LB$-$UB$ canonical form} 
whenever (i)--(iv) hold below.
\begin{enumerate}
\item ${\mathcal A}=\hbox{Mat}_{d+1}(\K)$.
\item $A$ is lower bidiagonal and $A^*$ is upper bidiagonal.
\item 
 $A_{i,i-1}=1 $ for $1 \leq i \leq d$.
\item $A_{00}=\theta_0$ and $A^*_{00}=\theta^*_0$.
\end{enumerate}
\end{definition}

\begin{lemma}
\label{thm:lbbasisshape}
Let $\Phi$ denote the Leonard system from 
(\ref{eq:ourstartingpt}).
Assume $\Phi$ is 
in $LB$-$UB$ canonical form, so that
 ${\mathcal A}=\hbox{Mat}_{d+1}(\K)$ by Definition
\ref{def:lbub}(i).
For $0 \leq i \leq d$ let $v_i$ denote the vector
in $\K^{d+1}$ which has ith coordinate  $1$
and all other coordinates 
0. Then the sequence $v_0, v_1, \ldots, v_d$ is a 
$\Phi$-$LB$-$UB$
basis for $\K^{d+1}$.
Let $\natural $ denote the  
$LB$-$UB$ canonical map for $\Phi$, from Definition
\ref{def:natcon}. Then $\natural $ is the identity map.
\end{lemma}
\noindent {\it Proof:}
Let $\theta_0, \theta_1, \ldots, \theta_d$ 
(resp. 
$\theta^*_0, \theta^*_1, \ldots, \theta^*_d$)
denote the eigenvalue
sequence (resp. dual eigenvalue sequence) for 
$\Phi$.
By Definition
\ref{def:lbub},
$A$ is lower bidiagonal with 
$A_{i,i-1}=1$ for $1 \leq i \leq d$.
By Corollary 
\ref{lem:basicdata}(i)
and since $A_{00}=\theta_0$ we find
 $A_{ii}=\theta_i$
for $0 \leq i \leq d$.
Apparently
$Av_i = \theta_iv_i+v_{i+1}$ for $0 \leq i \leq d-1$.
By Definition
\ref{def:lbub},
$A^*$  is upper
bidiagonal  with $A^*_{00}=\theta^*_0$.
Apparently $v_0$
 is an eigenvector for $A^*$ with eigenvalue $\theta^*_0$. Therefore 
$v_0 \in E^*_0V$.
Applying Lemma 
\ref{lem:lbubcheck}
 (with $V=\K^{d+1}$) 
we find 
 $v_0, v_1, \ldots, v_d$ is 
a $\Phi$-$LB$-$UB$
 basis for
$\K^{d+1}$. From
the construction 
each element in
 $\hbox{Mat}_{d+1}(\K)$ represents itself with respect to
 $v_0, v_1, \ldots, v_d$. Therefore $\natural $ is the identity map
 in view of Definition
\ref{def:natcon}.
\hfill $\Box $\\

\begin{theorem}
\label{thm:lbubshape}
Let $\Phi$ denote the Leonard
system 
from (\ref{eq:ourstartingpt})
and assume $\Phi$ is in $LB$-$UB$ canonical form.
Then $A=p^L$ and $A^*=p^U$, where $L,U$ are  from
Definition
\ref{def:LB}
and 
$p$
is the parameter array of $\Phi$.
\end{theorem}
\noindent {\it Proof:}
Let $\natural $ denote the $LB$-$UB$ canonical map for $\Phi$,
from Definition
\ref{def:natcon}. We assume $\Phi$ is in $LB$-$UB$ canonical form,
so $\natural $ is the identity map by Lemma
\ref{thm:lbbasisshape}.
Applying Lemma 
\ref{lem:lbubobv}
we find
$A=p^L$ and $A^*=p^U$.
\hfill $\Box $\\

\begin{corollary}
\label{cor:lbubiso}
Let $\Phi$ and $\Phi'$ denote Leonard systems over $\K$
which are in $LB$-$UB$ canonical form. Then the following
are equivalent: (i) $\Phi$ and $\Phi'$ are isomorphic; 
(ii) 
 $\Phi=\Phi'$.
\end{corollary}
\noindent {\it Proof:}
(i) $\Rightarrow $ (ii)
The Leonard systems
$\Phi, \Phi'$ have a common parameter array
which we denote by $p$. By Theorem 
\ref{thm:lbubshape} the Leonard pair associated with each of
$\Phi, \Phi'$ is equal to $p^L,p^U$. Apparently
$\Phi $ and $\Phi'$ are in the same associate class. By this
and since 
$\Phi, \Phi'$ are isomorphic we find $\Phi=\Phi'$ in view of 
Lemma
\ref{lem:nonisoprov}.
\\
(ii) $\Rightarrow $ (i) Clear.
\hfill $\Box $\\

\begin{definition}
\label{def:lbubconform}
\rm
Let $\Phi$ denote the 
Leonard system
from (\ref{eq:ourstartingpt}). By an {\it $LB$-$UB$ canonical form
for $\Phi$} we mean a Leonard system over $\K$ which is isomorphic
to $\Phi$ and which is in 
$LB$-$UB$ canonical form.
\end{definition}

\begin{theorem}
\label{def:lbconmap}
Let $\Phi$ denote the 
Leonard system
from (\ref{eq:ourstartingpt}). 
Then there exists a unique 
$LB$-$UB$ canonical form for $\Phi$.  
This form is $\Phi^\natural$, 
where $\natural $ denotes the
$LB$-$UB$ canonical map for $\Phi$ from Definition
\ref{def:natcon}.
\end{theorem}
\noindent {\it Proof:}
We first show $\Phi^\natural$ is an 
 $LB$-$UB$ canonical form for $\Phi$.
Since $\Phi$ is a Leonard system in $\mathcal A$
and since
$\natural : {\mathcal A}\rightarrow 
\hbox{Mat}_{d+1}(\K)$ is an isomorphism of $\K$-algebras,
we find
$\Phi^\natural$ is a Leonard system in
$\hbox{Mat}_{d+1}(\K)$ which is isomorphic to $\Phi$.
 We show 
$\Phi^\natural$ is in $LB$-$UB$ canonical  form.
To do this we show $\Phi^\natural$ satisfies
conditions (i)--(iv) of Definition
\ref{def:lbub}. Observe 
$\Phi^\natural$ satisfies
Definition
\ref{def:lbub}(i) since 
$\hbox{Mat}_{d+1}(\K)$ is the ambient algebra
of $\Phi^\natural$. 
Observe 
$\Phi^\natural$ satisfies
Definition
\ref{def:lbub}(ii)--(iv)
by
Definition \ref{def:LB} and
Lemma 
\ref{lem:lbubobv}.
We have now shown $\Phi^\natural$
satisfies 
Definition
\ref{def:lbub}(i)--(iv)
so $\Phi^\natural $  
is in 
$LB$-$UB$ canonical form. 
Apparently $\Phi^\natural$ is a Leonard system over $\K$
which is isomorphic to $\Phi$ and which is in $LB$-$UB$ canonical
form. Therefore
$\Phi^\natural$ is an $LB$-$UB$ canonical form for $\Phi$
by Definition
\ref{def:lbubconform}.
To finish the proof we
let  $\Phi'$ denote an $LB$-$UB$ canonical form for $\Phi$ and
show $\Phi'=\Phi^\natural$. Observe $\Phi', \Phi^\natural$
are isomorphic since they are both isomorphic to $\Phi$.
The Leonard systems $\Phi', \Phi^\natural $ are isomorphic 
and in $LB$-$UB$ canonical form
so $\Phi'=\Phi^\natural$ by Corollary
\ref{cor:lbubiso}.
\hfill $\Box $\\

\begin{corollary}
\label{cor:lbublstopa}
Consider the set
of Leonard systems over $\K$ which are in  
 $LB$-$UB$ canonical form.
We give a bijection from this set
 to the set of parameter arrays
over $\K$.
The bijection sends each Leonard system to its own parameter array.
\end{corollary}
\noindent {\it Proof:}
By 
the remark following
Definition
\ref{def:pals},
the map which sends a given Leonard system 
 to its parameter array induces a bijection
from the set of isomorphism
classes of
Leonard systems 
over $\K$ to the set of parameter
arrays over $\K$. By
Theorem \ref{def:lbconmap}
each of these isomorphism classes
contains a unique element which is in $LB$-$UB$ canonical
form. The result follows.
\hfill $\Box $\\

\section{The $LB$-$UB$ canonical form for Leonard pairs}

\noindent 
In this section we define and discuss the $LB$-$UB$ canonical form
for Leonard pairs. 
We begin with a comment.

\begin{lemma}
\label{lem:atmostonelb}
Let $A,A^*$ denote the Leonard pair from
Definition \ref{def:lpabs}. Then there exists
at most one Leonard system which is associated with
$A,A^*$ and which is in $LB$-$UB$ canonical form.
\end{lemma}
\noindent {\it Proof:}
Let $\Phi$ and $\Phi'$ denote
Leonard systems
which are  associated with $A,A^*$ and which are in
$LB$-$UB$ canonical form.
We show $\Phi=\Phi'$.
Since $\Phi, \Phi'$ are in the same associate class,
this will follow once we show
$\Phi, \Phi'$ have the same eigenvalue sequence and the same
dual eigenvalue sequence.
Observe by Theorem
\ref{thm:lbubshape} that the sequence of diagonal entries for
$A$ is the common eigenvalue sequence for 
$\Phi, \Phi'$.
Similary the sequence of diagonal entries for $A^*$ is
the common dual eigenvalue sequence for
$\Phi, \Phi'$.
Apparently $\Phi=\Phi'$.
\hfill $\Box $\\

\noindent Referring to the above lemma, we now
consider those Leonard pairs 
for which there exists an associated Leonard system
which is in 
$LB$-$UB$ canonical form.
In order to describe these we introduce the  
$LB$-$UB$ canonical form for Leonard pairs.

\begin{definition}
\label{lem:whenlb}
\rm
Let $A,A^*$ denote the Leonard pair from 
Definition \ref{def:lpabs}.
We say this pair is in 
{\it $LB$-$UB$ canonical form} 
whenever 
(i)--(iii) hold below.
\begin{enumerate}
\item ${\mathcal A}=\hbox{Mat}_{d+1}(\K)$.
\item $A$ is lower bidiagonal and $A^*$ is upper bidiagonal.
\item $A_{i,i-1}=1$ for $1 \leq i \leq d$.
\end{enumerate}
\end{definition}

\noindent We 
just defined the $LB$-$UB$ canonical form for Leonard pairs,
and in
 Definition
\ref{def:lbub}
we defined this form for Leonard systems.
We now compare these two versions. 
We will use the following definition. 

\begin{definition}
\label{def:desig}
\rm
Let $d$ denote a nonnegative integer and let 
$A,A^*$ denote a Leonard pair in 
$\hbox{Mat}_{d+1}(\K)$. We assume
$A$ is lower bidiagonal and $A^*$ is upper bidiagonal.
We make some comments and definitions.
(i) By Corollary  
\ref{lem:basicdata}(i)
the sequence $A_{00}, A_{11}, \ldots, A_{dd}$ is an eigenvalue
sequence for 
$A,A^*$. We call
this sequence
the {\it designated eigenvalue sequence} for $A,A^*$.
(ii)
By Corollary  
\ref{lem:basicdata}(iii)
 the sequence $A^*_{00}, A^*_{11}, \ldots, A^*_{dd}$ is a dual
 eigenvalue sequence for $A,A^*$.
We call this sequence
the {\it designated dual eigenvalue sequence} for
$A,A^*$.
(iii) By the {\it designated 
Leonard system} for $A,A^*$ we mean the
Leonard system
which 
is associated with $A,A^*$ and which has 
eigenvalue sequence 
 $A_{00}, A_{11}, \ldots, A_{dd}$
and dual eigenvalue sequence 
 $A^*_{00}, A^*_{11}, \ldots, A^*_{dd}$.
(iv) By the 
 {\it designated parameter array} for $A,A^*$ we mean the
 parameter array of the designated Leonard system for $A,A^*$.
\end{definition}

\begin{lemma}
\label{lem:lbubcompare1}
Let $A,A^*$ denote the Leonard pair from
Definition \ref{def:lpabs}.
Then the following are equivalent:
\begin{enumerate}
\item
$A,A^*$ is in $LB$-$UB$ canonical form.
\item There exists a Leonard system $\Phi$ which is associated with
$A,A^*$ and which is in 
$LB$-$UB$ canonical form.
\end{enumerate}
Suppose (i), (ii) hold. Then $\Phi$ is the designated Leonard
system of $A,A^*$.
\end{lemma}
\noindent {\it Proof:}
(i) $\Rightarrow $ (ii) Let $\Phi$ denote the designated Leonard
system for $A,A^*$, from Definition
\ref{def:desig}(iii). From the construction $\Phi$ is 
associated with $A,A^*$ and in $LB$-$UB$ canonical form.
\\
(ii) $\Rightarrow $ (i) 
Compare Definition 
\ref{def:lbub}
and Definition
\ref{lem:whenlb}.
\\
\noindent Now suppose (i), (ii) hold. Then $\Phi$ is
the designated Leonard system for $A,A^*$ by 
Lemma 
%Definition \ref{def:desig2}
\ref{lem:atmostonelb}
and the proof of
(i) $\Rightarrow $ (ii) above.
\hfill $\Box $\\

\begin{corollary}
\label{lem:lbubhowtell}
 We give a bijection from the set of
Leonard systems over $\K$ which are in $LB$-$UB$ canonical form,
to the set of Leonard pairs over $\K$ 
which are in $LB$-$UB$ canonical
form. The bijection sends each Leonard system to its associated
Leonard pair. The inverse bijection sends each
Leonard pair 
to its designated Leonard system.
\end{corollary}
\noindent {\it Proof:}
This is a reformulation of Lemma
\ref{lem:lbubcompare1}.
\hfill $\Box $\\

\begin{theorem}
\label{thm:lbubbijpa}
We give a bijection from the set of parameter arrays
over $\K$ to the set of Leonard pairs over $\K$ which
are  in $LB$-$UB$ canonical form. The bijection sends
each parameter array $p$ to the Leonard pair
$p^L, p^U$. The inverse bijection sends each
Leonard pair to its designated parameter array.
\end{theorem}
\noindent {\it Proof:}
Composing the inverse of the bijection from Corollary
\ref{cor:lbublstopa}, with the bijection from
Corollary
\ref{lem:lbubhowtell},
we obtain a bijection from
 the set of parameter arrays
over $\K$ to the set of Leonard pairs over $\K$ which
are  in $LB$-$UB$ canonical form. 
Let $p$ denote a parameter array over $\K$
and let $A,A^*$ denote the image of $p$ under this
bijection.
We show $A=p^L$ and $A^*=p^U$. 
By Corollary
\ref{cor:lbublstopa}
there exists a unique Leonard system 
over $\K$ which
is in $LB$-$UB$ canonical form and which has parameter array
$p$. Let us denote this system  by $\Phi$.
By the construction 
$A,A^*$
is associated with 
$\Phi$. 
Applying Theorem
\ref{thm:lbubshape} to $\Phi$ 
we find 
$A=p^L$ and $A^*=p^U$. 
To finish the proof we show $p$ is the designated parameter
array for $A,A^*$. We mentioned
$A,A^*$ is associated with
$\Phi$ 
 and $\Phi$ is in $LB$-$UB$ canonical form
so  
$\Phi$ is the designated Leonard system for $A,A^*$ by
Corollary
\ref{lem:lbubhowtell}.
We mentioned $p$ is the parameter array for $\Phi$ so
$p$ is the designated parameter array for $A,A^*$ by
Definition
\ref{def:desig}(iv).
\hfill $\Box $\\

\begin{definition}
\rm
Let $A,A^*$ denote the Leonard pair from 
Definition \ref{def:lpabs}. 
By an
{\it $LB$-$UB$ canonical form for $A,A^*$ }
we mean a Leonard pair over $\K$ which is isomorphic
to $A,A^*$ and which is in $LB$-$UB$ canonical form.
\end{definition}

\begin{theorem}
\label{thm:lbubcform4}
Let $A,A^*$ denote the Leonard pair from 
Definition \ref{def:lpabs}. 
We give a bijection from the set of
parameter arrays for $A,A^*$ to the
set of 
$LB$-$UB$ canonical forms for $A,A^*$.
This bijection sends each parameter array
$p$ to the pair $p^L, p^U$. (The parameter arrays for
$A,A^*$ are given in
Lemma
\ref{lem:palprel}.)
The inverse bijection sends each 
$LB$-$UB$ canonical form for $A,A^*$ to 
 its
designated parameter array.
\end{theorem}
\noindent {\it Proof:}
Let $B,B^*$ denote a Leonard pair over $\K$ which is
in $LB$-$UB$ canonical form. Let $p$ denote the designated
parameter array for $B,B^*$.
In view of Theorem 
\ref{thm:lbubbijpa} it suffices to show the following
are equivalent: 
(i) $A,A^*$ and $B,B^*$ are isomorphic;
(ii) $p$ is a parameter array for $A,A^*$.
These statements are equivalent by Lemma
\ref{lem:parrayforlp}.
\hfill $\Box $\\

\begin{corollary}
\label{cor:howmanylbubcf}
Let $A,A^*$ denote the Leonard pair from
Definition \ref{def:lpabs}. If $d\geq 1$ then
there exist exactly four $LB$-$UB$ canonical forms
for $A,A^*$. If $d=0$ there exists a unique
 $LB$-$UB$ canonical form
for $A,A^*$.
\end{corollary}
\noindent {\it Proof:}
Immediate from 
Theorem 
\ref{thm:lbubcform4}
and 
Corollary
\ref{cor:numberpa}.
\hfill $\Box $\\

\section{How to recognize a Leonard pair in $LB$-$UB$ canonical form} 

\noindent 
Let $d$ denote a nonnegative integer and let $A, A^*$ denote
matrices in $\hbox{Mat}_{d+1}(\K)$.
Let us assume 
$A$ is lower bidiagonal 
and $A^*$ is upper bidiagonal. 
We give a necessary and sufficient condition for
 $A,A^*$ to be a Leonard pair  which is in $LB$-$UB$ canonical
form.

\begin{theorem}
\label{thm:luc}
Let $d$ denote a nonnegative integer and let $A, A^*$ denote
matrices in $\hbox{Mat}_{d+1}(\K)$. Assume $A$ is lower bidiagonal
and $A^*$ is upper bidiagonal. Then the following (i), (ii) are equivalent.
\begin{enumerate}
\item The pair $A,A^*$ is a Leonard pair in $\hbox{Mat}_{d+1}(\K)$
which is in $LB$-$UB$ canonical form.
\item
There exists a parameter array
$(\theta_i, \theta^*_i, i=0..d;  \varphi_j, \phi_j, j=1..d)$
over $\K$ 
such that
\begin{eqnarray}
A_{ii} &=&\theta_i, \qquad \qquad 
A^*_{ii} =\theta^*_i \qquad \qquad (0 \leq i \leq d),
\label{eqn:entry1}
\\
A_{i,i-1} &=& 1, \qquad \qquad
A^*_{i-1,i} = \varphi_i \qquad \qquad (1 \leq i \leq d).
\label{eqn:entry2}
\end{eqnarray}
\end{enumerate}
Suppose (i), (ii) hold. Then the parameter array in (ii) above 
is uniquely determined by $A,A^*$.
This parameter array is 
the designated parameter array for 
$A,A^*$ in the sense of 
Definition
\ref{def:desig}.
\end{theorem}
\noindent {\it Proof:}
This is a reformulation of Theorem
\ref{thm:lbubbijpa}.
\hfill $\Box $\\

\section{Leonard pairs $A,A^*$ with $A$ lower bidiagonal and
$A^*$ upper bidiagonal}

\noindent 
Let $d$ denote a nonnegative integer and let $A, A^*$ denote
matrices in $\hbox{Mat}_{d+1}(\K)$. Let us assume 
$A$ is lower bidiagonal 
and $A^*$ is upper  bidiagonal. 
We give a necessary and sufficient condition
for $A,A^*$ to be  a Leonard pair.

\begin{theorem}
\label{thm:lug}
Let $d$ denote a nonnegative integer and let $A, A^*$ denote
matrices in $\hbox{Mat}_{d+1}(\K)$. Assume $A$  
lower bidiagonal and $A^*$ is upper bidiagonal.
Then the following (i), (ii) are equivalent.
\begin{enumerate}
\item
The pair $A,A^*$ is a Leonard pair in 
$\hbox{Mat}_{d+1}(\K)$.
\item
There exists a parameter array
$(\theta_i, \theta^*_i, i=0..d;  \varphi_j, \phi_j, j=1..d)$
over $\K$ 
such that
\begin{eqnarray}
&&A_{ii} =\theta_i,
\qquad \qquad 
A^*_{ii} =\theta^*_i \qquad \qquad (0 \leq i \leq d),
\label{eq:comb1}
\\
&&\qquad A_{i,i-1}A^*_{i-1,i} = \varphi_i \qquad \qquad (1 \leq i \leq d).
\label{eq:comb2}
\end{eqnarray}
\end{enumerate}
Suppose (i), (ii) hold. Then the parameter array in (ii) above 
is uniquely determined by $A,A^*$.
This parameter array is the designated parameter array
for $A,A^*$ in the sense of 
Definition
\ref{def:desig}.
\end{theorem}
\noindent {\it Proof:}
\noindent 
(i) $\Rightarrow $ (ii)
By Corollary 
\ref{lem:basicdata}(ii)
we have 
 $A_{i,i-1} \not=0$ for $1 \leq i \leq d$.
Let $S$ denote the diagonal matrix in
 $\hbox{Mat}_{d+1}(\K)$
which has diagonal entries 
$S_{ii}=A_{10}A_{21}\cdots A_{i,i-1}$
for $0 \leq i \leq d$.
Each of $S_{00}, S_{11}, \ldots, S_{dd}$ 
is nonzero 
 so $S^{-1}$ exists.
Let
$
\sigma :  
 \hbox{Mat}_{d+1}(\K) \rightarrow 
 \hbox{Mat}_{d+1}(\K)
$
denote the isomorphism of $\K$-algebras which satisfies
$X^\sigma = 
 S^{-1}XS$
for all $X \in  
 \hbox{Mat}_{d+1}(\K)$.
From the construction
$A^\sigma $ (resp. $A^{*\sigma} $)
is lower bidiagonal (resp. upper bidiagonal)
with entries
\begin{eqnarray}
&&A^\sigma_{ii}=A_{ii}, \qquad \qquad A^{*\sigma}_{ii} = A^*_{ii} 
\qquad \qquad 
(0 \leq i \leq d),
\label{eq:bshape1}
\\
&&A^\sigma_{i,i-1}=1, \qquad \qquad 
A^{*\sigma}_{i-1,i}=A_{i,i-1}A^*_{i-1,i}
\qquad \qquad  (1 \leq i \leq d).
\label{eq:bshape2}
\end{eqnarray}
Apparently 
$A^\sigma,A^{*\sigma}$ is a Leonard pair in
 $\hbox{Mat}_{d+1}(\K)$ 
 which is 
in $LB$-$UB$ canonical form.
Applying Theorem 
\ref{thm:luc} to this pair 
we find there exists a parameter
array
$(\theta_i, \theta^*_i, i=0..d;  \varphi_j, \phi_j, j=1..d)$
over $\K$ 
such that both
$A^\sigma_{ii} =\theta_i
$, $A^{*\sigma}_{ii} =\theta^*_i$ for $0 \leq i \leq d$
and 
$A^{*\sigma}_{i-1,i} = \varphi_i$ for  $1 \leq i \leq d$.
Combining  these facts with
(\ref{eq:bshape1}), 
(\ref{eq:bshape2}) 
we find this parameter array satisfies
(\ref{eq:comb1}),
(\ref{eq:comb2}).
\\
(ii) $\Rightarrow $ (i)
For $1 \leq i \leq d$ we have $A_{i,i-1} \not=0$
by 
(\ref{eq:comb2}) and since $\varphi_i\not=0$.
Let  
$
\sigma :  
 \hbox{Mat}_{d+1}(\K) \rightarrow 
 \hbox{Mat}_{d+1}(\K)
$
denote the isomorphism of $\K$-algebras from the proof of
(i) $\Rightarrow $ (ii) above.
We routinely
 find both
$A^\sigma_{ii} =\theta_i$,
$A^{*\sigma}_{ii} =\theta^*_i$ for $0 \leq i \leq d$ and
both $A^\sigma_{i,i-1} = 1$,
$A^{*\sigma}_{i-1,i} = \varphi_i$ for $1 \leq i \leq d$.
Apparently
$A^{\sigma},A^{*\sigma}$ satisfies 
Theorem
\ref{thm:luc}(ii). Applying that theorem  to this pair
we find
$A^\sigma, A^{*\sigma}$ is a Leonard pair in 
 $\hbox{Mat}_{d+1}(\K)$
which is in
$LB$-$UB$ canonical form. In particular
$ A^\sigma, A^{*\sigma}$ is a Leonard pair
in  $\hbox{Mat}_{d+1}(\K)$. By this and since
$\sigma$ is an
isomorphism we find
$A,A^*$ is a Leonard pair in
 $\hbox{Mat}_{d+1}(\K)$.

\medskip
\noindent Suppose (i), (ii) hold  above.
Let $p$ denote a parameter array which satisfies  
(ii) above.
We show $p$ is the designated parameter array
for $A,A^*$.
We first show $p$ is a parameter array for $A,A^*$.
Observe 
$p$ is a parameter array for
$A^\sigma, A^{*\sigma}$ 
by Theorem
\ref{thm:luc}
and the proof of 
(ii) $\Rightarrow $ (i) above.
Also $A,A^*$ is isomorphic to   
 $A^\sigma, A^{*\sigma}$ so $p$ is a parameter array for
$A,A^*$. Observe $p$ is the 
 designated parameter array for $A,A^*$  
by Definition
\ref{def:desig}.
\hfill $\Box $\\

\section{Examples of Leonard pairs $A,A^*$ with $A$ lower bidiagonal
and $A^*$ upper bidiagonal}

\begin{example}
\label{ex:lbublp}
\rm
Let $d$ denote a nonnegative integer. Let $A$ and $A^*$ denote
the following matrices in $\hbox{Mat}_{d+1}(\K)$.
\beast
A=\left(
\begin{array}{ c c c c c c}
d &   &      &      &   &{\bf 0} \\
-1 & d-2  &     &      &   &  \\
  & -2  &  \cdot    &   &   & \\
  &   & \cdot     & \cdot  &    & \\
  &   &           &  \cdot & \cdot &  \\
{\bf 0} &   &   &   & -d & -d  
\end{array}
\right),
\qquad
A^*=\left(
\begin{array}{ c c c c c c}
d & 2d  &      &      &   &{\bf 0} \\
 & d-2  & 2d-2    &      &   &  \\
  &   &  \cdot    & \cdot  &   & \\
  &   &      & \cdot  &  \cdot  & \\
  &   &           &   & \cdot & 2 \\
{\bf 0} &   &   &   &  & -d  
\end{array}
\right) .
\eeast
To avoid degenerate situations, 
we assume the characteristic of $\K$ is zero or an odd prime
greater than $d$. 
Then the pair $A,A^*$ is a Leonard pair 
in 
$\hbox{Mat}_{d+1}(\K)$.
The corresponding designated 
parameter array
from
Definition \ref{def:desig} 
 is the parameter array given
in Example \ref{ex:pa1}.
\end{example}
\noindent {\it Proof:}
Let 
$(\theta_i, \theta^*_i, i=0..d;  \varphi_j, \phi_j, j=1..d)$
denote the parameter array from Example
\ref{ex:pa1}.
We routinely find this parameter array satisfies 
Theorem 
\ref{thm:lug}(ii);
applying that theorem
we find $A,A^*$ is a Leonard pair in 
$\hbox{Mat}_{d+1}(\K)$.
The parameter array
$(\theta_i, \theta^*_i, i=0..d;  \varphi_j, \phi_j, j=1..d)$
is the designated parameter array of $A,A^*$ by the last line of
Theorem 
\ref{thm:lug}.
\hfill $\Box $\\

\begin{example}
\label{ex:lbublp2}
\rm
Let $d,
q, s, s^*, r_1, r_2$ be as in Example
\ref{ex:pa2}.
Let $A$ and $A^*$ denote the following matrices in
 $\hbox{Mat}_{d+1}(\K)$.
The matrix $A$ is lower bidiagonal
with entries
\beast
A_{ii} &=& q^{-i} + sq^{i+1} \qquad \qquad (0 \leq i \leq d),
\\
A_{i,i-1} &=& (1-q^{-i})(1-r_1q^i) \qquad \quad (1 \leq i \leq d).
\eeast
The matrix $A^*$ 
is upper bidiagonal
 with entries
\beast
A^*_{ii} &=& q^{-i} + s^*q^{i+1} \qquad \qquad (0 \leq i \leq d),
\\
A^*_{i-1,i} &=& (q^{-d}-q^{1-i})(1-r_2q^i)
\qquad \quad (1 \leq i \leq d).
\eeast
Then the pair $A,A^*$ is a Leonard pair in
$\hbox{Mat}_{d+1}(\K)$.
The corresponding designated parameter array 
from Definition \ref{def:desig} 
is the parameter array given in Example \ref{ex:pa2}.
\end{example}
\noindent {\it Proof:}
Let 
$(\theta_i, \theta^*_i, i=0..d;  \varphi_j, \phi_j, j=1..d)$
denote the parameter array from
Example \ref{ex:pa2}. 
We routinely find this array satisfies 
Theorem 
\ref{thm:lug}(ii);
applying that
theorem 
we find $A,A^*$ is a Leonard pair in 
$\hbox{Mat}_{d+1}(\K)$.
The parameter array
$(\theta_i, \theta^*_i, i=0..d;  \varphi_j, \phi_j, j=1..d)$
is the designated parameter array for $A,A^*$ by the last line
of Theorem
\ref{thm:lug}.
\hfill $\Box $\\

\section{The $TD$-$D$ canonical form; preliminaries}

\medskip
\noindent 
We now turn our attention to 
the $TD$-$D$
 canonical form.
We begin with some comments.

\begin{lemma}
 \cite[Lemma 5.1]{LS24}
\label{lem:tddbasis}
Let $\Phi$ denote the Leonard system from 
(\ref{eq:ourstartingpt}) and let 
$V$ denote an irreducible $\mathcal A$-module.
Let $\eta_0$ denote a nonzero vector in $E_0V$. 
Then the sequence
\begin{eqnarray}
E^*_0\eta_0, E^*_1\eta_0, \ldots, E^*_d\eta_0
\label{eq:dst}
\end{eqnarray}
is a basis for $V$.
\end{lemma}

\begin{definition} 
\label{def:tddbasis}
\rm
Let $\Phi$ denote the Leonard system from 
(\ref{eq:ourstartingpt}) and let 
$V$ denote an irreducible  $\mathcal A$-module.
By a {\it $\Phi$-$TD$-$D$ basis}
for $V$ we mean 
a sequence of the form
(\ref{eq:dst}), where $\eta_0$ denotes a nonzero vector
in $E_0V$.
\end{definition}

\medskip
\noindent 
The concept of a $\Phi$-$TD$-$D$  basis will play an important
role in what follows. 
Therefore we examine it carefully. In each of the next two
lemmas we give a
characterization of this type of basis.

\begin{lemma}
\label{lem:eggechar}
Let $\Phi$ denote the Leonard system from 
(\ref{eq:ourstartingpt}) and let 
$V$ denote an irreducible $\mathcal A$-module.
Let $v_0, v_1, \ldots, v_d$ denote a sequence of vectors in
$V$, not all 0. Then  this sequence is a 
$\Phi$-$TD$-$D$ basis for $V$ 
 if and only if both (i), (ii) hold below.
\begin{enumerate}
\item $v_i \in E^*_iV$ for $0 \leq i \leq d$.
\item $\sum_{i=0}^dv_i\in E_0V$.
\end{enumerate}
\end{lemma}
\noindent {\it Proof:}
To prove the lemma in one direction,
 assume 
$v_0, v_1, \ldots, v_d$ is a
$\Phi$-$TD$-$D$ basis for $V$.
By Definition
\ref{def:tddbasis}
there exists a nonzero $\eta_0 \in E_0V$ such
that $v_i = E^*_i\eta_0$ for $0 \leq i \leq d$.
Apparently $v_i \in E^*_iV$ for $0 \leq i \leq d$ so (i) holds.
Let $I$ denote the identity element of $\mathcal A$ and
observe $I=\sum_{i=0}^d E^*_i$.
 Applying this to $\eta_0$ we find
$\eta_0=\sum_{i=0}^d v_i $ and (ii) follows.
We have now proved the lemma in one direction. To prove the lemma
in the other direction, assume $v_0, v_1, \ldots, v_d$ satisfy
(i), (ii) above.
We define
$\eta_0=\sum_{i=0}^d v_i$ and observe $\eta_0 \in E_0V$.
Using (i) 
we find $E^*_iv_j=\delta_{ij}v_j$ for
$0 \leq i,j\leq d$; it follows
$v_i = E^*_i\eta_0$ for $0 \leq i \leq d$.
Observe $\eta_0 \not=0$  since at least one
of $v_0, v_1, \ldots, v_d$ is nonzero.
Now 
$v_0, v_1, \ldots, v_d$ is 
a 
$\Phi$-$TD$-$D$  basis for $V$ by Definition
\ref{def:tddbasis}.
\hfill $\Box $\\

\noindent  We recall some notation.
Let $d$ denote a nonnegative
integer and let $B$ denote a matrix in 
$\hbox{Mat}_{d+1}(\K)$. Let $\alpha $ denote a scalar in $\K$.
Then $B$ is said to have {\it constant row  sum $\alpha$} whenever
$B_{i0}+B_{i1}+\cdots + B_{id}=\alpha$ for $0 \leq i \leq d$.

\begin{lemma}
\label{lem:rowsum}
Let $\Phi$ denote the Leonard system from 
(\ref{eq:ourstartingpt}).
Let $\theta_0, \theta_1, \ldots, \theta_d$ 
(resp. 
$\theta^*_0, \theta^*_1, \ldots$, $\theta^*_d$) 
denote  
the eigenvalue sequence (resp. dual eigenvalue sequence) of $\Phi$.
Let $V$ denote an irreducible $\mathcal A$-module and
let $v_0, v_1, \ldots, v_d$ denote a basis
for $V$.
Let $B$ (resp. $B^*$) denote the matrix in 
$\hbox{Mat}_{d+1}(\K)$
which represents
$A$ (resp. $A^*)$ with respect to 
this basis.
Then 
 $v_0, v_1, \ldots, v_d$ is a $\Phi$-$TD$-$D$ basis
 for $V$ if and only if (i), (ii) hold below.
\begin{enumerate}
\item $B$ has constant row sum $\theta_0$.
\item $B^*=\hbox{diag}(\theta^*_0, \theta^*_1, \ldots, \theta^*_d)$.
\end{enumerate}
\end{lemma}
\noindent {\it Proof:}
Observe $A \sum_{j=0}^d v_j = \sum_{i=0}^d v_i(B_{i0}+B_{i1}+\cdots B_{id})$.
Recall $E_0V$ is the eigenspace for $A$ and eigenvalue $\theta_0$.
Apparently 
$B$ has constant row sum $\theta_0$ if and only if
$\sum_{i=0}^d v_i \in E_0V$.
Recall that for $0 \leq i \leq d$,
$E^*_iV$ is the eigenspace for $A^*$ and eigenvalue $\theta^*_i$.
Apparently 
$B^*=\hbox{diag}(\theta^*_0, \theta^*_1, \ldots, \theta^*_d)$
if and only if $v_i \in E^*_iV$ for $0 \leq i \leq d$. 
The result follows in view of Lemma
\ref{lem:eggechar}.
\hfill $\Box $\\

\section{The $TD$-$D$ canonical map}

Let $\Phi$ denote the Leonard system from  
(\ref{eq:ourstartingpt}). In this section we use
$\Phi$ to define
a certain isomorphism 
$\flat : {\mathcal A}\rightarrow 
\hbox{Mat}_{d+1}(\K)$. We  call $\flat$
the {\it $TD$-$D$ canonical map} for $\Phi$.
We describe the entries of $A^\flat$ and $A^{*\flat}$.

\begin{definition}
\label{def:flatcon}
\rm
Let $\Phi$ denote the Leonard system from 
(\ref{eq:ourstartingpt}). We define
a map $\flat : {\mathcal A}\rightarrow 
\hbox{Mat}_{d+1}(\K)$ as follows.
Let 
$V$ denote an irreducible $\mathcal A$-module.
For all $X \in {\mathcal  A}$ we let $X^\flat $ denote
the matrix 
in 
$\hbox{Mat}_{d+1}(\K)$
which represents $X$ with respect to 
a $\Phi$-$TD$-$D$ basis for $V$.
We observe $\flat : {\mathcal A} \rightarrow 
\hbox{Mat}_{d+1}(\K)$ is an isomorphism of
$\K$-algebras. We call $\flat $ the
{\it $TD$-$D$ canonical map} for $\Phi$.
\end{definition}

\noindent 
Referring to Definition
\ref{def:flatcon},
we now describe $A^\flat $ and $A^{*\flat}$.
We begin with a comment.

\begin{lemma}
\label{lem:firstcom}
Let $\Phi$ denote the Leonard system from 
(\ref{eq:ourstartingpt}). 
Let $\theta_0, \theta_1, \ldots, \theta_d$ 
(resp. 
$\theta^*_0, \theta^*_1, \ldots$, $\theta^*_d$) 
denote  
the eigenvalue sequence (resp. dual eigenvalue sequence) of $\Phi$.
Let $\flat $ denote the $TD$-$D$
canonical map for $\Phi$, from Definition
\ref{def:flatcon}.
Then (i), (ii) hold below.
\begin{enumerate}
\item $A^\flat $ has constant row sum $\theta_0$.
\item $A^{*\flat}=\hbox{diag}(\theta^*_0,\theta^*_1, \ldots, \theta^*_d)$.
\end{enumerate}
\end{lemma}
\noindent {\it Proof:}
Combine 
Lemma 
\ref{lem:rowsum} and
Definition
\ref{def:flatcon}.
\hfill $\Box $\\

\noindent 
Referring to Definition
\ref{def:flatcon}, we now  
describe  $A^\flat $
and $A^{*\flat}$ from another point of view. 
 We use the following notation.

\begin{definition}
\label{def:TD}
\rm
Consider the set of all parameter arrays
over $\K$. 
We define two functions on this set.
 We call these functions 
 $T$ and $D$.
Let 
$p=(\theta_i, \theta^*_i, i=0..d;  \varphi_j, \phi_j, j=1..d)$
denote a parameter array over $\K$.
The image $p^T$ is the tridiagonal matrix in
$\hbox{Mat}_{d+1}(\K)$ 
which has 
the following entries. The diagonal entries are
\beast
p^T_{ii} = \theta_i + \frac{\varphi_i}{\theta^*_i-\theta^*_{i-1}}
 + \frac{\varphi_{i+1}}{\theta^*_i-\theta^*_{i+1}}
\eeast
for $0 \leq i \leq d$, where we recall 
$\varphi_0=0$, $\varphi_{d+1}=0$ and where
$\theta^*_{-1}, \theta^*_{d+1}$ denote indeterminates.
The superdiagonal and subdiagonal entries are
\beast
p^{T}_{i-1,i}= 
\varphi_i \frac{\prod_{h=0}^{i-2}(\theta^*_{i-1}-\theta^*_h)
}
{\prod_{h=0}^{i-1}(\theta^*_{i}-\theta^*_h)
},
\qquad \qquad 
p^{T}_{i,i-1}=
 \phi_{i} \frac{\prod_{h=i+1}^d (\theta^*_i-\theta^*_h)
}
{\prod_{h=i}^d(\theta^*_{i-1}-\theta^*_h)
}
\eeast
for $1 \leq i \leq d$. 
The image $p^D$ is the following  matrix in 
$\hbox{Mat}_{d+1}(\K)$.
\beast
p^D=\hbox{diag}(\theta^*_0, \theta^*_1, \ldots, \theta^*_d).
\eeast
\end{definition}

\begin{theorem}
\label{thm:quote}
Let $\Phi$ denote the 
Leonard system
from (\ref{eq:ourstartingpt}).
Let $\flat $ denote the $TD$-$D$
canonical map for $\Phi$, from Definition
\ref{def:flatcon}.
Then
$A^\flat=p^T$ and $A^{*\flat}=p^D$, where 
$p$ denotes the parameter array for $\Phi$.
\end{theorem}
\noindent {\it Proof:}
Observe 
$A^{*\flat}=p^D$ by 
Lemma
\ref{lem:firstcom}(ii).
We have 
$A^\flat=p^T$  by  
\cite[Theorem 11.2]{LS24}.
\hfill $\Box $\\

\noindent We finish this section with an observation. 

\begin{corollary}  
\label{cor:aipbipci}
Let 
$p=(\theta_i, \theta^*_i, i=0..d;  \varphi_j, \phi_j, j=1..d)$
denote a parameter array over $\K$. Then
the matrix $p^{T}$ has constant row sum $\theta_0$.
\end{corollary}
\noindent {\it Proof:}
By the remark after Definition
\ref{def:pals} there exists a Leonard system
$\Phi$ over 
 $\K$ which has
parameter array $p$. 
For notational convenience let us assume $\Phi$ is 
the Leonard system
(\ref{eq:ourstartingpt}).
Let $\flat $ denote the $TD$-$D$ canonical map for
$\Phi$, from
Definition
\ref{def:flatcon}.
Then $A^\flat $ has constant row sum $\theta_0$ by
Lemma
\ref{lem:firstcom} and
$A^\flat = p^T$ by Theorem
\ref{thm:quote} so
 $p^T$ has constant row
sum $\theta_0$.
\hfill $\Box $\\

\section{The $TD$-$D$ canonical form for Leonard systems} 

\noindent 
In this section we introduce the 
$TD$-$D$ canonical form for Leonard systems.
We 
define what it means for a given Leonard system
 to be in $TD$-$D$ canonical
form. We describe the 
  Leonard systems which are in 
$TD$-$D$ canonical form.
 We  show 
every Leonard system is isomorphic to 
a unique Leonard system 
which is in $TD$-$D$ canonical form.

\begin{definition}
\label{def:tddcon}
\rm
Let $\Phi$ denote the Leonard system from 
(\ref{eq:ourstartingpt}).
Let $\theta_0, \theta_1, \ldots, \theta_d$ 
(resp. 
$\theta^*_0, \theta^*_1, \ldots$, $\theta^*_d$) 
denote  
the eigenvalue sequence (resp. dual eigenvalue sequence) of $\Phi$.
We say  $\Phi$ is in 
{\it $TD$-$D$ canonical form}
whenever (i)--(iii) hold below.
\begin{enumerate}
\item ${\mathcal A}=\hbox{Mat}_{d+1}(\K)$.
\item $A$ is tridiagonal and $A^*$ is diagonal.
\item $A$ has constant row sum $\theta_0$ and $A^*_{00}=\theta^*_0$.
\end{enumerate}
\end{definition}

\begin{lemma}
\label{thm:tddbasisshape}
Let $\Phi$ denote the Leonard system from 
(\ref{eq:ourstartingpt}).
Assume $\Phi$ is 
in $TD$-$D$ canonical form, so that
 ${\mathcal A}=\hbox{Mat}_{d+1}(\K)$ by Definition
\ref{def:tddcon}(i).
For $0 \leq i \leq d$ let $v_i$ denote the vector
in $\K^{d+1}$ which has ith coordinate  $1$
and all other coordinates 
0. Then the sequence $v_0, v_1, \ldots, v_d$ is a 
$\Phi$-$TD$-$D$
basis for $\K^{d+1}$.
Let $\flat $ denote the  
$TD$-$D$ canonical map for $\Phi$, from Definition
\ref{def:flatcon}. Then $\flat $ is the identity map.
\end{lemma}
\noindent {\it Proof:}
Observe $v_0, v_1, \ldots, v_d$ is a basis for
$\K^{d+1}$, and that with respect to this basis each element
of 
 $\hbox{Mat}_{d+1}(\K)$ represents itself. 
Let $\theta^*_0, \theta^*_1, \ldots, \theta^*_d$
denote the dual eigenvalue sequence for $\Phi$.
By 
Corollary
\ref{cor:tdeig}(ii) and since $A^*_{00}=\theta^*_0$
we find $A^*=\hbox{diag}(\theta^*_0, \theta^*_1,\ldots, \theta^*_d)$.
Applying 
Lemma
\ref{lem:rowsum} (with $V=\K^{d+1}$), 
we find 
 $v_0, v_1, \ldots, v_d$ is 
a $\Phi$-$TD$-$D$
 basis for
$\K^{d+1}$.
We mentioned each element in
 $\hbox{Mat}_{d+1}(\K)$ represents itself with respect to
 $v_0, v_1, \ldots, v_d$, so $\flat $ is the identity map
 in view of Definition
\ref{def:flatcon}.
\hfill $\Box $\\

\begin{theorem}
\label{thm:tdshape}
Let $\Phi$ denote the Leonard
system 
from (\ref{eq:ourstartingpt}),
and assume $\Phi$ is in $TD$-$D$ canonical form.
Then $A=p^T$ and $A^*=p^D$, where $T,D$ are from
Definition
\ref{def:TD}
and 
$p$
is the parameter array for $\Phi$.
\end{theorem}
\noindent {\it Proof:}
Let $\flat $ denote the $TD$-$D$ canonical map for $\Phi$,
from Definition
\ref{def:flatcon}. We assume $\Phi$ is in $TD$-$D$ canonical form,
so $\flat $ is the identity map by Lemma
\ref{thm:tddbasisshape}.
Applying 
Theorem \ref{thm:quote} we find
$A=p^T$ and $A^*=p^D$.
\hfill $\Box $\\

\begin{corollary}
\label{cor:tddiso}
Let $\Phi$ and $\Phi'$ denote Leonard systems over $\K$
which are in $TD$-$D$ canonical form. Then the following
are equivalent: (i) $\Phi$ and $\Phi'$ are isomorphic; 
(ii) 
 $\Phi=\Phi'$.
\end{corollary}
\noindent {\it Proof:}
(i) $\Rightarrow $ (ii)
The Leonard systems
$\Phi, \Phi'$ have a common parameter array
which we denote by $p$. By Theorem 
\ref{thm:tdshape} the Leonard pair associated with each of
$\Phi, \Phi'$ is equal to $p^T,p^D$. Apparently
$\Phi $ and $\Phi'$ are in the same associate class. By this
and since 
$\Phi, \Phi'$ are isomorphic we find $\Phi=\Phi'$ in view of 
Lemma
\ref{lem:nonisoprov}.
\\
(ii) $\Rightarrow $ (i) Clear.
\hfill $\Box $\\

\begin{definition}
\label{def:tddconform}
\rm
Let $\Phi$ denote the 
Leonard system
from (\ref{eq:ourstartingpt}). By a {\it $TD$-$D$ canonical form
for $\Phi$} we mean a Leonard system over $\K$ which is isomorphic
to $\Phi$ and which is in 
$TD$-$D$ canonical form.
\end{definition}

\begin{theorem}
\label{lem:stmap}
\label{def:tddconmap}
Let $\Phi$ denote the 
Leonard system
from (\ref{eq:ourstartingpt}). 
Then there exists a unique 
$TD$-$D$ canonical form for $\Phi$.  
This is $\Phi^\flat$, 
where $\flat $ denotes the
$TD$-$D$ canonical map for $\Phi$ from Definition
\ref{def:flatcon}.
\end{theorem}
\noindent {\it Proof:}
We first show $\Phi^\flat$ is a 
 $TD$-$D$ canonical form for $\Phi$.
Since $\Phi$ is a Leonard system in $\mathcal A$
and since
$\flat : {\mathcal A}\rightarrow 
\hbox{Mat}_{d+1}(\K)$ is an isomorphism of $\K$-algebras,
we find
$\Phi^\flat$ is a Leonard system in
$\hbox{Mat}_{d+1}(\K)$ which is isomorphic to $\Phi$.
 We show 
$\Phi^\flat$ is in $TD$-$D$ canonical  form.
To do this we show $\Phi^\flat$ satisfies
conditions (i)--(iii) of Definition
\ref{def:tddcon}. Observe 
$\Phi^\flat$ satisfies
Definition
\ref{def:tddcon}(i) since 
$\hbox{Mat}_{d+1}(\K)$ is the ambient algebra
of $\Phi^\flat$. 
Observe 
$\Phi^\flat$ satisfies
Definition
\ref{def:tddcon}(ii),(iii)
by
Lemma \ref{lem:firstcom}
and 
Theorem \ref{thm:quote}.
We have now shown $\Phi^\flat$
satisfies 
Definition
\ref{def:tddcon}(i)--(iii)
so $\Phi^\flat $  
is in 
$TD$-$D$ canonical form. 
Apparently $\Phi^\flat$ is a Leonard system over $\K$
which is isomorphic to $\Phi$ and which is in $TD$-$D$ canonical
form. Therefore
$\Phi^\flat$ is a $TD$-$D$ canonical form for $\Phi$
by Definition
\ref{def:tddconform}.
To finish the proof we
let  $\Phi'$ denote a $TD$-$D$ canonical form for $\Phi$ and
show $\Phi'=\Phi^\flat$. Observe $\Phi', \Phi^\flat$
are isomorphic since they are both isomorphic to $\Phi$.
The Leonard systems $\Phi', \Phi^\flat $are isomorphic 
and in $TD$-$D$ canonical form
so $\Phi'=\Phi^\flat$ by Corollary
\ref{cor:tddiso}.
\hfill $\Box $\\

\begin{corollary}
\label{cor:lstopa} Consider the set
of Leonard systems over $\K$ which are in  
 $TD$-$D$ canonical form.
We give a bijection from this set
 to the set of parameter arrays
over $\K$.
The bijection sends each Leonard system to its own parameter array.
\end{corollary}
\noindent {\it Proof:}
By 
the remark following
Definition
\ref{def:pals},
the map which sends a given Leonard system 
 to its parameter array induces a bijection
from the set of isomorphism
classes of
Leonard systems 
over $\K$ to the set of parameter
arrays over $\K$. By
Theorem \ref{def:tddconmap}
each of these isomorphism classes
contains a unique element which is in $TD$-$D$ canonical
form. The result follows.
\hfill $\Box $\\

\section{The $TD$-$D$ canonical form for Leonard pairs}

\noindent 
In this section we define and discuss the $TD$-$D$ canonical form
for Leonard pairs. 
We begin with a comment.

\begin{lemma}
\label{lem:atmostone}
Let $A,A^*$ denote the Leonard pair from
Definition \ref{def:lpabs}. Then there exists
at most one Leonard system which is associated with
$A,A^*$ and which is in $TD$-$D$ canonical form.
\end{lemma}
\noindent {\it Proof:}
Let $\Phi$ and $\Phi'$ denote
Leonard systems
which are  associated with $A,A^*$ and which are in
$TD$-$D$ canonical form.
We show $\Phi=\Phi'$.
Let 
$\theta_0, \theta_1, \ldots, \theta_d$ 
(resp.
$\theta'_0, \theta'_1, \ldots, \theta'_d$)
denote the eigenvalue sequence for $\Phi$ (resp. $\Phi'$.)
Let 
$\theta^*_0, \theta^*_1, \ldots, \theta^*_d$ 
(resp.
$\theta^{*\prime}_0, \theta^{*\prime}_1, \ldots, \theta^{*\prime}_d$)
denote the dual eigenvalue sequence for $\Phi$ (resp. $\Phi'$.)
Observe $\Phi, \Phi'$ are in the same associate class so
$\Phi'$ is one of $\Phi, \Phi^\downarrow, \Phi^\Downarrow, 
\Phi^{\downarrow\Downarrow}$.
Therefore $\theta'_i = \theta_i$ for $0 \leq i \leq d$
or 
 $\theta'_i = \theta_{d-i}$ for $0 \leq i \leq d$.
Also 
$\theta^{*\prime}_i = \theta^*_i$ for $0 \leq i \leq d$
or 
 $\theta^{*\prime}_i = \theta^*_{d-i}$ for $0 \leq i \leq d$.
To show $\Phi=\Phi'$ it suffices to show
 $\theta_i = \theta'_i$ and 
$\theta^{*}_i = \theta^{*\prime}_i$ for $0 \leq i \leq d$.
Each of $\theta_0, \theta'_0$ is equal to the common row
sums of $A$ so $\theta_0=\theta'_0$. Apparently
 $\theta_i = \theta'_i$ for $0 \leq i\leq d$.
Each of 
 $\theta^*_0, \theta^{*\prime}_0$ is equal to
$A^*_{00}$ so 
 $\theta^*_0=\theta^{*\prime}_0$.
Apparently
$\theta^{*}_i = \theta^{*\prime}_i$ for $0 \leq i \leq d$.
We conclude $\Phi=\Phi'$.
\hfill $\Box $\\

\noindent Referring to the above lemma, we now
consider those Leonard pairs 
for which there exists an associated Leonard system
which is in 
$TD$-$D$ canonical form.
In order to describe these we introduce the  
$TD$-$D$ canonical form for Leonard pairs.

\begin{definition}
\label{def:lptdd}
\rm
Let $A,A^*$ denote the Leonard pair from  
Definition \ref{def:lpabs} and let
$\theta_0, \theta_1, \ldots, \theta_d$
denote an  eigenvalue sequence for this pair.
We say $A,A^*$ is in {\it $TD$-$D$ canonical form}
whenever (i)--(iii) hold below.
\begin{enumerate}
\item ${\mathcal A}=\hbox{Mat}_{d+1}(\K)$.
\item $A$ is tridiagonal and $A^*$ is diagonal.
\item $A$ has constant row sum and this sum is $\theta_0$ or $\theta_d$.
\end{enumerate}
\end{definition}

\noindent We 
just defined the $TD$-$D$ canonical form for Leonard pairs,
and in
 Definition
\ref{def:tddcon}
we defined this form for Leonard systems.
We now compare these two versions. 
We will use the following definition. 

\begin{definition}
\label{def:desig2}
\rm
Let $A,A^*$ denote the Leonard pair from
Definition \ref{def:lpabs}, and assume 
this pair is in 
 $TD$-$D$ 
canonical form. We make several comments and definitions.
(i) By Definition 
\ref{def:lptdd}(iii) and 
Definition
\ref{def:evseqlp}, there exists a unique 
eigenvalue sequence 
$\theta_0, \theta_1,\ldots, \theta_d$
for $A,A^*$
such that $A$ has constant row sum $\theta_0$.
We call this the {\it designated eigenvalue sequence} for $A,A^*$. 
(ii)
By Corollary 
\ref{cor:tdeig}(ii)
the sequence
 $A^*_{00}, A^*_{11}, \ldots, A^*_{dd}$ is a dual eigenvalue
 sequence for $A,A^*$. We call this the 
{\it designated dual eigenvalue sequence} for
$A,A^*$.
(iii) By the {\it designated Leonard system} for $A,A^*$ we mean the
Leonard system which is associated with $A,A^*$ and which has
 eigenvalue sequence 
$\theta_0, \theta_1,\ldots, \theta_d$
and dual eigenvalue sequence
 $A^*_{00}, A^*_{11}, \ldots, A^*_{dd}$.
(iv)  By the {\it designated parameter array} for $A,A^*$ we mean the
parameter array of the designated Leonard system for $A,A^*$.
\end{definition}

\begin{lemma}
\label{lem:compare1}
Let $A,A^*$ denote the Leonard pair from
Definition \ref{def:lpabs}.
Then the following are equivalent:
\begin{enumerate}
\item
$A,A^*$ is in $TD$-$D$ canonical form.
\item There exists a Leonard system $\Phi$ which is associated with
$A,A^*$ and which is in 
$TD$-$D$ canonical form.
\end{enumerate}
Suppose (i), (ii) hold. Then $\Phi$ is the designated Leonard
system of $A,A^*$.
\end{lemma}
\noindent {\it Proof:}
(i) $\Rightarrow $ (ii) Let $\Phi$ denote the designated Leonard
system for $A,A^*$, from Definition
\ref{def:desig2}(iii). From the construction $\Phi$ is 
associated with $A,A^*$ and in $TD$-$D$ canonical form.
\\
(ii) $\Rightarrow $ (i)
Compare Definition
\ref{def:tddcon}
and 
 Definition \ref{def:lptdd}. 
\\
\noindent Now suppose (i), (ii) hold. Then $\Phi$ is
the designated Leonard system for $A,A^*$ by 
Lemma 
\ref{lem:atmostone}
and the proof of
(i) $\Rightarrow $ (ii) above.
\hfill $\Box $\\

\begin{corollary}
\label{lem:howtelltdd}
 We give a bijection from the set of
Leonard systems over $\K$ which are in $TD$-$D$ canonical form,
to the set of Leonard pairs over $\K$ 
which are in $TD$-$D$ canonical
form. The bijection sends each Leonard system to its associated
Leonard pair. The inverse bijection sends each
Leonard pair 
to its designated Leonard system.
\end{corollary}
\noindent {\it Proof:}
This is a reformulation of Lemma
\ref{lem:compare1}.
\hfill $\Box $\\

\begin{theorem}
\label{thm:bijTDpa}
We give a bijection from the set of parameter arrays
over $\K$ to the set of Leonard pairs over $\K$ which
are  in $TD$-$D$ canonical form. The bijection sends
each parameter array $p$ to the Leonard pair
$p^T, p^D$. The inverse bijection sends each
Leonard pair to its designated parameter array.
\end{theorem}
\noindent {\it Proof:}
Composing the inverse of the bijection from Corollary
\ref{cor:lstopa}, with the bijection from
Corollary
\ref{lem:howtelltdd},
we obtain a bijection from
 the set of parameter arrays
over $\K$ to the set of Leonard pairs over $\K$ which
are  in $TD$-$D$ canonical form. 
Let $p$ denote a parameter array over $\K$
and let $A,A^*$ denote the image of $p$ under this
bijection.
We show $A=p^T$ and $A^*=p^D$. 
By Corollary
\ref{cor:lstopa}
there exists a unique Leonard system 
over $\K$ which
is in $TD$-$D$ canonical form and which has parameter array
$p$. Let us denote this system  by $\Phi$.
By the construction 
$A,A^*$
is associated with 
$\Phi$. 
Applying Theorem
\ref{thm:tdshape} to $\Phi$ 
we find 
$A=p^T$ and $A^*=p^D$. 
To finish the proof we show $p$ is the designated parameter
array for $A,A^*$. We mentioned
$A,A^*$ is associated with
$\Phi$ 
 and $\Phi$ is in $TD$-$D$ canonical form
so  
$\Phi$ is the designated Leonard system for $A,A^*$ by
Corollary
\ref{lem:howtelltdd}.
We mentioned $p$ is the parameter array for $\Phi$ so
$p$ is the designated parameter array for $A,A^*$ by
Definition
\ref{def:desig2}(iv).
\hfill $\Box $\\

\begin{definition}
\rm
Let $A,A^*$ denote the Leonard pair from 
Definition \ref{def:lpabs}. 
By a {\it $TD$-$D$ canonical form for $A,A^*$} 
we mean a Leonard pair over $\K$ which is isomorphic
to $A,A^*$ and which is in $TD$-$D$ canonical form.
\end{definition}

\begin{theorem}
\label{thm:tddcform4}
Let $A,A^*$ denote the Leonard pair from 
Definition \ref{def:lpabs}. 
We give a bijection from the set of
parameter arrays for $A,A^*$ to the
set of 
$TD$-$D$ canonical forms for $A,A^*$.
This bijection sends each parameter array
$p$ to the pair $p^T, p^D$. (The parameter arrays for
$A,A^*$ are given in
Lemma
\ref{lem:palprel}.)
The inverse bijection sends each 
$TD$-$D$ canonical form for $A,A^*$ to 
 its
designated parameter array.
\end{theorem}
\noindent {\it Proof:}
Let $B,B^*$ denote a Leonard pair over $\K$ which is
in $TD$-$D$ canonical form. Let $p$ denote the designated
parameter array for $B,B^*$.
In view of Theorem 
\ref{thm:bijTDpa} it suffices to show the following
are equivalent: 
(i) $A,A^*$ and $B,B^*$ are isomorphic;
(ii) $p$ is a parameter array for $A,A^*$.
These statements are equivalent by Lemma
\ref{lem:parrayforlp}.
\hfill $\Box $\\

\begin{corollary}
\label{cor:howmanytddcf}
Let $A,A^*$ denote the Leonard pair from
Definition \ref{def:lpabs}. If $d\geq 1$ then
there exist exactly four $TD$-$D$ canonical forms
for $A,A^*$. If $d=0$ then there exists a unique
 $TD$-$D$ canonical form
for $A,A^*$.
\end{corollary}
\noindent {\it Proof:}
Immediate from 
Theorem 
\ref{thm:tddcform4}
and 
Corollary
\ref{cor:numberpa}.
\hfill $\Box $\\

\section{How to recognize a Leonard pair in $TD$-$D$ canonical form}

\noindent 
Let $d$ denote a nonnegative integer and let $A, A^*$ denote
matrices in $\hbox{Mat}_{d+1}(\K)$.
Let us assume 
$A$ is tridiagonal 
and $A^*$ is diagonal. 
We give a necessary and sufficient condition for  
$A,A^*$ to be a Leonard pair which is in $TD$-$D$ canonical
form. We present two versions of our result.

\begin{theorem}
\label{thm:tdstand}
Let $d$ denote a nonnegative integer and let $A, A^*$ denote
matrices in $\hbox{Mat}_{d+1}(\K)$. Assume $A$  is 
tridiagonal and $A^*$ is diagonal.
Then the following (i), (ii) are equivalent.
\begin{enumerate}
\item
The pair $A,A^*$ is a Leonard pair in
 $\hbox{Mat}_{d+1}(\K)$ which is in
$TD$-$D$ canonical form.
\item
There exists
a parameter array 
$(\theta_i, \theta^*_i, i=0..d;  \varphi_j, \phi_j, j=1..d)$
over $\K$
such that
\beast
A_{ii} &=&
\theta_i + \frac{\varphi_i}{\theta^*_i-\theta^*_{i-1}}
+
 \frac{\varphi_{i+1}}{\theta^*_i-\theta^*_{i+1}}
\qquad 
(0 \leq i \leq d),
\\
A_{i-1,i}&=& 
\varphi_i \frac{\prod_{h=0}^{i-2}(\theta^*_{i-1}-\theta^*_h)
}
{\prod_{h=0}^{i-1}(\theta^*_{i}-\theta^*_h)
} \qquad \qquad (1 \leq i \leq d),
\\
A_{i,i-1}&=&
 \phi_{i} \frac{\prod_{h=i+1}^d (\theta^*_i-\theta^*_h)
}
{\prod_{h=i}^d(\theta^*_{i-1}-\theta^*_h)
}
\qquad \qquad (1 \leq i \leq d),
\\
A^*_{ii} &=&\theta^*_i \qquad \qquad (0 \leq i \leq d).
\eeast
\end{enumerate}
Suppose (i), (ii) hold. Then the parameter array in (ii) above 
is uniquely determined by $A,A^*$.
This parameter array is the designated parameter array
for $A,A^*$ in the sense of 
Definition
\ref{def:desig2}.
\end{theorem}
\noindent {\it Proof:}
This is a reformulation of
Theorem
\ref{thm:bijTDpa}.
\hfill $\Box $\\

\begin{theorem}
\label{thm:tdstand2}
Let $d$ denote a nonnegative integer and let $A, A^*$ denote
matrices in $\hbox{Mat}_{d+1}(\K)$. Assume $A$   is
tridiagonal and $A^*$ is diagonal.
Then the following (i), (ii) are equivalent.
\begin{enumerate}
\item
The pair $A,A^*$ is a Leonard pair in
 $\hbox{Mat}_{d+1}(\K)$ which is in
$TD$-$D$ canonical form.
\item  
There exists
a parameter array 
$(\theta_i, \theta^*_i, i=0..d;  \varphi_j, \phi_j, j=1..d)$
over $\K$
such that $A$ has constant row sum $\theta_0$ and
\beast
A_{i-1,i}&=& 
\varphi_i \frac{\prod_{h=0}^{i-2}(\theta^*_{i-1}-\theta^*_h)
}
{\prod_{h=0}^{i-1}(\theta^*_{i}-\theta^*_h)
} \qquad \qquad (1 \leq i \leq d),
\\
A_{i,i-1}&=&
 \phi_{i} \frac{\prod_{h=i+1}^d (\theta^*_i-\theta^*_h)
}
{\prod_{h=i}^d(\theta^*_{i-1}-\theta^*_h)
}
\qquad \qquad (1 \leq i \leq d),
\\
A^*_{ii} &=&\theta^*_i \qquad \qquad (0 \leq i \leq d).
\eeast
\end{enumerate}
Suppose (i), (ii) hold. Then the parameter array in
(ii) above is uniquely determined by $A,A^*$. This parameter
array is the designated parameter array for $A,A^*$ in
the sense of 
Definition
\ref{def:desig2}.
\end{theorem}
\noindent {\it Proof:}
Combine Theorem 
\ref{thm:tdstand} and Corollary
\ref{cor:aipbipci}.
\hfill $\Box $\\

\section{Examples of Leonard pairs in $TD$-$D$ canonical form}

\noindent In this section we give a few examples of Leonard
pairs which are in $TD$-$D$ canonical form.

\begin{example}
\label{ex:kraw} 
\rm
Let $d$ denote a nonnegative integer. 
Let $A$ and $A^*$ denote the following matrices in $\hbox{Mat}_{d+1}(\K)$.

\beast
A = 
\left(
\begin{array}{ c c c c c c}
0 & d  &      &      &   &{\bf 0} \\
1 & 0  &  d-1   &      &   &  \\
  & 2  &  \cdot    & \cdot  &   & \\
  &   & \cdot     & \cdot  & \cdot   & \\
  &   &           &  \cdot & \cdot & 1 \\
{\bf 0} &   &   &   & d & 0  
\end{array}
\right),
\qquad A^*= \hbox{diag}(d, d-2, d-4, \ldots, -d).
\eeast
To avoid degenerate situations,
we assume the characteristic of $\K$ is zero or an odd prime
greater than $d$.
Then 
the pair $A,A^*$ is a Leonard
pair in 
 $\hbox{Mat}_{d+1}(\K)$
which is in
$TD$-$D$ canonical form. The corresponding   
designated
 parameter array 
from Definition \ref{def:desig2}
is the 
parameter array
given in Example
\ref{ex:pa1}.
\end{example}
\noindent {\it Proof:}
Let 
$(\theta_i, \theta^*_i, i=0..d;  \varphi_j, \phi_j, j=1..d)$
denote the parameter array from Example
\ref{ex:pa1}.
We routinely verify this parameter
array  
 satisfies 
Theorem \ref{thm:tdstand2}(ii); applying that theorem
we find $A,A^*$ is a Leonard pair in 
 $\hbox{Mat}_{d+1}(\K)$
which is in
$TD$-$D$ canonical form.
The parameter array
$(\theta_i, \theta^*_i, i=0..d;  \varphi_j, \phi_j, j=1..d)$
is the designated parameter array for $A,A^*$ by the last line
of 
Theorem \ref{thm:tdstand2}.
\hfill $\Box $\\

\begin{example}
\label{ex:qrac3}
\rm
Let 
 $d,q, s, s^*, r_1, r_2$ be as in
Example
 \ref{ex:pa2}.
 Let $A$ and $A^*$
  denote the following matrices in $\hbox{Mat}_{d+1}(\K)$.
The matrix $A$ is tridiagonal
 with entries
\beast
A_{01} &=& \frac{(1-q^{-d})(1-r_1q)(1-r_2q)}
{1-s^*q^2},
\\
A_{i-1,i} &=& \frac{(1-q^{i-d-1})(1-s^*q^i)(1-r_1q^i)(1-r_2q^i)}
{(1-s^*q^{2i-1})(1-s^*q^{2i})} \qquad \quad (2 \leq i \leq d),
\\
A_{i,i-1} &=& \frac{(1-q^i)(1-s^*q^{i+d+1})(r_1-s^*q^i)(r_2-s^*q^i)}
{s^*q^d(1-s^*q^{2i})(1-s^*q^{2i+1})} \qquad (1 \leq i \leq d-1),
\\
A_{d,d-1} &=& \frac{(1-q^d)(r_1-s^*q^d)(r_2-s^*q^d)}
{s^*q^d(1-s^*q^{2d})}
\eeast
and  constant row sum $1+sq$. 
The matrix $A^*$ is diagonal 
 with entries
\beast
A^*_{ii} &=& q^{-i}+s^*q^{i+1}\qquad (0 \leq i \leq d).
\eeast
Then the pair $A,A^*$ is a Leonard pair in
 $\hbox{Mat}_{d+1}(\K)$ which is in
$TD$-$D$ canonical form.
The corresponding designated  parameter array 
from Definition \ref{def:desig2}
 is the  parameter array given in Example \ref{ex:pa2}.
\end{example}
\noindent {\it Proof:}
Let
$(\theta_i, \theta^*_i, i=0..d;  \varphi_j, \phi_j, j=1..d)$
denote the parameter array from Example \ref{ex:pa2}.
We routinely verify this parameter
array  
 satisfies 
Theorem \ref{thm:tdstand2}(ii); applying that theorem
we find $A,A^*$ is a Leonard pair in 
 $\hbox{Mat}_{d+1}(\K)$ 
which is in $TD$-$D$ canonical form.
The parameter array
$(\theta_i, \theta^*_i, i=0..d;  \varphi_j, \phi_j, j=1..d)$
is the designated parameter array for $A,A^*$ by the last line
of 
Theorem \ref{thm:tdstand2}.
\hfill $\Box $\\

\section{Leonard pairs $A,A^*$ with $A$ tridiagonal and $A^*$ diagonal}

\noindent 
Let $d$ denote a nonnegative integer and let $A, A^*$ denote
matrices in $\hbox{Mat}_{d+1}(\K)$. Let us assume 
$A$ is tridiagonal 
and $A^*$ is diagonal. 
We give  a necessary and sufficient condition for
 $A,A^*$ to be a Leonard pair. 

\begin{theorem}
\label{thm:tdcrit}
Let $d$ denote a nonnegative integer and let $A, A^*$ denote
matrices in $\hbox{Mat}_{d+1}(\K)$. Assume $A$ is
tridiagonal and $A^*$ is diagonal.
Then the following (i), (ii) are equivalent.
\begin{enumerate}
\item
The pair $A,A^*$ is a Leonard pair 
in 
 $\hbox{Mat}_{d+1}(\K)$.
\item
There exists
a parameter array 
$(\theta_i, \theta^*_i, i=0..d;  \varphi_j, \phi_j, j=1..d)$
over $\K$
such that
\begin{eqnarray}
A_{ii} &=&
\theta_i + \frac{\varphi_i}{\theta^*_i-\theta^*_{i-1}}
+
 \frac{\varphi_{i+1}}{\theta^*_i-\theta^*_{i+1}}
\qquad \qquad  
(0 \leq i \leq d), 
\label{eq:aiiform}
\\
A_{i,i-1}A_{i-1,i}&=&
\varphi_i\phi_i \frac{\prod_{h=0}^{i-2}(\theta^*_{i-1}-\theta^*_h)
}
{\prod_{h=0}^{i-1}(\theta^*_{i}-\theta^*_h)
}
\, \frac{\prod_{h=i+1}^d (\theta^*_i-\theta^*_h)
}
{\prod_{h=i}^d(\theta^*_{i-1}-\theta^*_h)
}
\qquad  (1 \leq i \leq d),
\label{eq:crossprod}
\\
A^*_{ii} &=&\theta^*_i \qquad \qquad (0 \leq i \leq d).
\label{eq:ths}
\end{eqnarray}
\end{enumerate}
Suppose 
 (i), (ii) hold and let $R$ denote the set
 of parameter arrays which satisfy 
(ii) above. Then $R$ consists of the parameter arrays
$(\theta_i, \theta^*_i, i=0..d;  \varphi_j, \phi_j, j=1..d)$
for $A,A^*$ which 
satisfy $\theta^*_i=A^*_{ii}$ for $0 \leq i \leq d$.
If 
$(\theta_i, \theta^*_i, i=0..d;  \varphi_j, \phi_j, j=1..d)$
is in $R$
then so is
$(\theta_{d-i}, \theta^*_i, i=0..d;  \phi_j, \varphi_j, j=1..d)$
 and $R$ contains no further elements.  
\end{theorem}
\noindent {\it Proof:}
(i) $\Rightarrow $ (ii)
We assume $A$ is tridiagonal and $A^*$ is diagonal
so 
 $A^*_{00}, A^*_{11}, \ldots, A^*_{dd}$ 
is a dual eigenvalue sequence for
$A,A^*$ by Corollary
\ref{cor:tdeig}(ii).
For notational convenience we set $\theta^*_i=A^*_{ii} $ for
$0\leq i \leq d$. 
By Definition 
\ref{def:evseqlp} there exists a
Leonard system 
$\Phi$ 
which is associated with $A,A^*$ and
which has dual eigenvalue sequence
$\theta^*_0, \theta^*_1, \ldots, \theta^*_d$.
Let 
$\theta_0, \theta_1, \ldots, \theta_d$ denote
the eigenvalue sequence for $\Phi$.
Let $\varphi_1, \varphi_2, \ldots, \varphi_d$ (resp. 
 $\phi_1, \phi_2, \ldots, \phi_d$) denote the first
 (resp. second) split sequence for $\Phi$.
We abbreviate
$p=(\theta_i, \theta^*_i, i=0..d;  \varphi_j, \phi_j, j=1..d)$
and observe $p$ 
is the parameter array for $\Phi$.
%by
%Definition
%\ref{def:pals}.
We show $p$ satisfies the conditions of (ii) above. 
Observe $p$  is over $\K$ since the 
Leonard pair $A,A^*$ is over $\K$.
We show $p$ satisfies 
(\ref{eq:aiiform})--(\ref{eq:ths}).
Let $\flat$
denote the $TD$-$D$ canonical map for $\Phi$. We recall
$A^\flat=p^T$ and $A^{*\flat}=p^D$ by 
Theorem \ref{thm:quote}.
Since 
$\hbox{Mat}_{d+1}(\K)$ is the ambient algebra
of $\Phi$ the domain 
of 
 $\flat$ is 
equal to $\hbox{Mat}_{d+1}(\K)$.
Since the range of 
 $\flat$ is 
equal to $\hbox{Mat}_{d+1}(\K)$ as well,
there exists an invertible
matrix $S \in 
\hbox{Mat}_{d+1}(\K)$
such that $X^\flat=SXS^{-1}$ for all $X \in 
\hbox{Mat}_{d+1}(\K)$. Observe $A^{*\flat}=A^*$ so
$SA^*=A^*S$. The matrix $A^*$ is diagonal 
with diagonal entries mutually distinct so
$S$ is diagonal. From this and since $A^\flat=SAS^{-1}$
we find
$
A^\flat_{ii}=
A_{ii}
$
for $0 \leq i \leq d$ and
$
A^\flat_{i,i-1}A^\flat_{i-1,i} 
=
A_{i,i-1}A_{i-1,i} 
$
for $1 \leq i \leq d$.
By these comments 
the parameter array $p$
satisfies 
(\ref{eq:aiiform}) and 
(\ref{eq:crossprod}). From the construction $p$
satisfies (\ref{eq:ths}).
\\
\noindent
(ii) $\Rightarrow $ (i)
Let 
$p:=(\theta_i, \theta^*_i, i=0..d;  \varphi_j, \phi_j, j=1..d)$
denote a parameter array over $\K$ which satisfies
(\ref{eq:aiiform})--(\ref{eq:ths}).
Let $\Phi$ denote a Leonard system over $\K$ 
which has parameter  array $p$.
Recall $\Phi$ is only determined up to isomorphism; replacing
$\Phi$ with an isomorphic Leonard system if necessary
we may assume  
$\Phi$ is in $TD$-$D$ canonical form by
Theorem \ref{lem:stmap}.
Let $B,B^*$ denote the Leonard pair associated with $\Phi$.
Then $B=p^T$ and $B^*=p^D$
by Theorem
\ref{thm:tdshape}.
Apparently
$B^*=A^*$; moreover
$B_{ii}=A_{ii}$ for $0 \leq i \leq d$ and
$B_{i,i-1}B_{i-1,i}=A_{i,i-1}A_{i-1,i}$ for $1 \leq i \leq d$.
Let $S$ denote the diagonal matrix in
$\hbox{Mat}_{d+1}(\K)$
which has  diagonal entries
$S_{ii}=\prod_{h=1}^i A_{i,i-1}/B_{i,i-1}$ for $0 \leq i \leq d$.
We observe  $S_{ii}\not=0$ for $0 \leq i \leq d$ so
$S^{-1}$ exists.
Let $\sigma :
\hbox{Mat}_{d+1}(\K)
\rightarrow \hbox{Mat}_{d+1}(\K)$
denote the isomorphism of $\K$-algebras which satisfies
$X^\sigma = SXS^{-1}$ for all $X \in 
\hbox{Mat}_{d+1}(\K)$. From our above comments we find
$B^\sigma=A$ and $
B^{*\sigma}=A^*$.
By this and since $B,B^*$ is a Leonard pair in
$\hbox{Mat}_{d+1}(\K)$ we find
$A,A^*$ is a Leonard pair
in $\hbox{Mat}_{d+1}(\K)$.
\\
\noindent Suppose (i), (ii) hold.
Let $R'$ denote the set of parameter
arrays
for $A,A^*$ which have dual eigenvalue sequence 
$A^*_{00}, A^*_{11}, \ldots, A^*_{dd}$. From Lemma
\ref{lem:palprel}
we find  
that if 
 $(\theta_i, \theta^*_i, i=0..d;  \varphi_j, \phi_j, j=1..d)$
is in $R'$ then so is
$(\theta_{d-i}, \theta^*_i, i=0..d;  \phi_j, \varphi_j, j=1..d)$
and $R'$ contains no further elements.
We now  show $R=R'$.
From the proof of
(i) $\Rightarrow $ (ii) above we find $R' \subseteq R$.
We show 
$R \subseteq R'$. Let 
 $(\theta_i, \theta^*_i, i=0..d;  \varphi_j, \phi_j, j=1..d)$
denote a parameter array in $R$.
 By the proof of
(ii) $\Rightarrow $ (i)
above we find this array is for $A,A^*$ in the sense of
Definition \ref{def:palp}.
By (\ref{eq:ths}) we find
$\theta^*_i=A^*_{ii}$ for $0 \leq i\leq d$. 
Apparently 
 $(\theta_i, \theta^*_i, i=0..d;  \varphi_j, \phi_j, j=1..d)$
is contained in 
$R'$ and it follows
$R \subseteq R'$. We have now shown $R=R'$ and the proof is complete.   
\hfill $\Box $\\

\section{How to compute the parameter arrays which satisfy Theorem 
\ref{thm:tdcrit}(ii) 
}

\noindent Let $d$ denote a positive integer and
let $A,A^*$ denote a Leonard pair in
$\hbox{Mat}_{d+1}(\K)$. Let us assume 
$A$ is tridiagonal and $A^*$ is diagonal. Suppose we wish
to verify that $A,A^*$ is a Leonard pair.
In order to do this it suffices to display 
a parameter array
$
(\theta_i, \theta^*_i, i=0..d;  \varphi_j, \phi_j, j=1..d)
$
which satisfies
Theorem
\ref{thm:tdcrit}(ii). 
We give a method for obtaining this array from the entries
of $A$ and $A^*$.
Our method is summarized as follows.
From 
(\ref{eq:ths}) we find $\theta^*_i=A^*_{ii}$ for $0 \leq i \leq d$.
To obtain the rest of the array
we proceed in two steps: (i) we obtain $\theta_0, \theta_d$
as the roots of a certain quadratic polynomial whose coefficients
are rational expressions involving $A_{00}, A_{11}, A_{dd}, A_{10}A_{01}$ and 
$\theta^*_0, \theta^*_1, \ldots, \theta^*_d$;
(ii) we obtain $\theta_i$ $(1 \leq i \leq d-1)$ and
$\varphi_i, \phi_i$ $(1 \leq i \leq d)$ as rational expressions involving 
$\theta_0, \theta_d, A_{00}, A_{dd}$ and  
$\theta^*_0, \theta^*_1, \ldots, \theta^*_d$.
For convenience  we discuss step (ii) before step (i).
To prepare for step (ii) we give a lemma. 

\begin{lemma}
\label{lem:moreeq}
Let $d$ denote a positive integer and let 
$
(\theta_i, \theta^*_i, i=0..d;  \varphi_j, \phi_j, j=1..d)
$
denote a parameter array over $\K$.
For notational convenience we define  
\begin{eqnarray}
\vartheta_i:= 
\sum_{h=0}^{i-1}
{{\theta^*_h-\theta^*_{d-h}}\over{\theta^*_0-\theta^*_d}} 
\qquad \qquad 
(0 \leq i \leq d).
\label{eq:vatheta}
\end{eqnarray}
Then  (i)--(iii) hold below.
\begin{enumerate}
\item 
$\displaystyle{
\theta_i =\theta_0 +
\frac{\varphi_i-\phi_d \vartheta_i}{\theta^*_{i-1}-\theta^*_d}
\qquad \qquad (1 \leq i \leq d).
}
$
\item 
$\displaystyle{
\theta_i =\theta_d +
\frac{\varphi_{i+1}-\phi_1 \vartheta_{i+1}}{\theta^*_{i+1}-\theta^*_0}
\qquad \qquad (0 \leq i \leq d-1).
}
$
\item 
$\displaystyle{
\frac{\varphi_{i+1}-\phi_1 \vartheta_{i+1}}{\theta^*_{i+1}-\theta^*_0}
=
\frac{\varphi_i-\phi_d \vartheta_i}{\theta^*_{i-1}-\theta^*_d}
+\theta_0-\theta_d \qquad \qquad (1 \leq i \leq d-1).
}
$
\end{enumerate}
\end{lemma}
\noindent {\it Proof:}
(i)
Let the integer $i$ be given.
Evaluating 
Corollary
\ref{cor:d4pa}(ii) using
Corollary
\ref{cor:d4pa}(i) 
we find
$\varphi_i=\phi_d \vartheta_i +
(\theta_i-\theta_0)(\theta^*_{i-1}-\theta^*_d)$. Solving
this equation for $\theta_i$ we get the result.
\\
\noindent (ii) Similar to the proof of (i) above,
except use Theorem
\ref{thm:classls}(iii) instead of 
Corollary
\ref{cor:d4pa}(ii).
\\
\noindent (iii) Combine (i), (ii) above.
\hfill $\Box $\\

\begin{theorem}
\label{lem:step2}
Let $d$ denote a positive  integer and let
 $A,A^*$ denote a Leonard pair in
$\hbox{Mat}_{d+1}(\K)$. Assume 
$A$ is tridiagonal and $A^*$ is diagonal.
Let
$(\theta_i, \theta^*_i, i=0..d;  \varphi_j, \phi_j, j=1..d)
$
denote a parameter array which satisfies
Theorem 
\ref{thm:tdcrit}(ii). 
Then
$\theta_i$ $(1 \leq i \leq d-1)$ and
$\varphi_i, \phi_i$ $(1 \leq i \leq d)$ 
are obtained from 
 $\theta_0, \theta_d, A_{00}, A_{dd}$ and
$\theta^*_0, \theta^*_1, \ldots, \theta^*_d$ as follows.
\begin{enumerate}
\item 
To obtain $\varphi_1, \varphi_d, \phi_1, \phi_d$ use
\begin{eqnarray}
&&
\varphi_1 = (A_{00}-\theta_0)(\theta^*_0-\theta^*_1),
\qquad \qquad 
\varphi_d = (A_{dd}-\theta_d)(\theta^*_d-\theta^*_{d-1}),
\label{eq:ends}
\\
&&
\phi_1 = (A_{00}-\theta_d)(\theta^*_0-\theta^*_1),
\qquad \qquad 
\phi_d = (A_{dd}-\theta_0)(\theta^*_d-\theta^*_{d-1}).
\label{eq:ends2}
\end{eqnarray}
\item 
To obtain
 $ \varphi_2,\varphi_3, \ldots, 
\varphi_{d-1}$ recursively apply 
Lemma
\ref{lem:moreeq}(iii).
\item To obtain $\theta_1, \theta_2, \ldots, \theta_{d-1}$ use
Lemma
\ref{lem:moreeq}(i) or Lemma \ref{lem:moreeq}(ii).
\item To obtain $\phi_2, \phi_3, \ldots, \phi_{d-1}$ use
Theorem 
\ref{thm:classls}(iv).  
\end{enumerate}
\end{theorem}
\noindent {\it Proof:}
(i)
To obtain 
the equation on the left (resp. right) in (\ref{eq:ends})
set $i=0$ (resp. $i=d$)  in 
(\ref{eq:aiiform}) and rearrange terms.
Line (\ref{eq:ends2}) is just
 (\ref{eq:ends}) 
with the original parameter array
replaced by the parameter array 
$(\theta_{d-i}, \theta^*_i, i=0..d;  \phi_j, \varphi_j, j=1..d)$.
\\
\noindent (ii)--(iv) Clear.
\hfill $\Box $\\

\begin{theorem}
\label{def:eps}
With reference to 
Theorem
\ref{lem:step2}, 
 the scalars $\theta_0, \theta_d$
 are the roots of the quadratic polynomial
\begin{eqnarray}
(\lambda-A_{00})(\lambda-\alpha/\varepsilon)-
A_{10}A_{01}/\varepsilon,
\label{eq:quad}
\end{eqnarray}
where  $\varepsilon, \alpha $ are defined as follows.
If $d=1$ then $\varepsilon =1$ and $\alpha = A_{11}$.
If $d\geq 2$ then 
\begin{eqnarray}
\varepsilon = 
 \frac{(\theta^*_1-\theta^*_d)
(\theta^*_1-\theta^*_{d-1}) \cdots 
(\theta^*_1-\theta^*_{2})} 
{(\theta^*_0-\theta^*_d)
(\theta^*_0-\theta^*_{d-1}) \cdots 
(\theta^*_0-\theta^*_{2})}
\label{eq:epsdef1}
\end{eqnarray}
and
\begin{eqnarray}
\alpha &=& A_{11} \frac{
\theta^*_1-\theta^*_2}
{\theta^*_0-\theta^*_2}
\;-
\;
A_{00}
\frac{
\theta^*_1-\theta^*_d
}
{
\theta^*_0-\theta^*_2
}
\,
\frac{\theta^*_0-\theta^*_1}{\theta^*_0-\theta^*_d}
\;+
\;
A_{dd}
\frac{
\theta^*_{d-1}-\theta^*_d}
{\theta^*_0-\theta^*_2}
\,\frac{\theta^*_0-\theta^*_1}{\theta^*_0-\theta^*_d}.
\label{eq:alphadef}
\end{eqnarray}
\end{theorem}
\noindent {\it Proof:}
First suppose $d=1$. Then
$\theta_0, \theta_d$ are the roots of the characteristic
polynomial of $A$ and this polynomial is $(\lambda - A_{00})(\lambda -A_{11})
-A_{10}A_{01}$. Next suppose $d\geq 2$.
We claim the scalar $\varepsilon$ from
(\ref{eq:epsdef1}) satisfies 
\begin{eqnarray}
 \varepsilon = 1 - \frac{\theta^*_0-\theta^*_1}{\theta^*_0-\theta^*_2}
\;
\frac{
\theta^*_0+\theta^*_1-\theta^*_{d-1}-\theta^*_d}{\theta^*_0-\theta^*_d}.
\label{eq:epsdef2}
\end{eqnarray}
 To obtain
(\ref{eq:epsdef2}) we 
recall by Corollary 
\ref{cor:aipbipci} 
that $p^T$ has constant row sum $\theta_0$, where
$
p=
(\theta_i, \theta^*_i, i=0..d;  \varphi_j, \phi_j, j=1..d)
$.
Considering row 1 of $p^T$ we find
$p^T_{10}+p^T_{11}+p^T_{12}=\theta_0$.
We evaluate the left-hand side of this equation
using Definition
\ref{def:TD}.
In the resulting equation 
% consider the equation of Corollary 
%\ref{cor:aipbipci} at $i=1$. In this equation we
we eliminate
$\varphi_1, \varphi_2$ using Theorem
\ref{thm:classls}(iii)   
and we simplify the result using Corollary
\ref{cor:d4pa}(i).
Line (\ref{eq:epsdef2}) follows and our claim is proved. 
To show $\theta_0, \theta_d$ are the roots of 
(\ref{eq:quad}) 
we show
both
\begin{eqnarray}
\theta_0+\theta_d &=& A_{00}+ \alpha/\varepsilon,
\label{eq:sum0d}
\\
\theta_0\theta_d&=&A_{00} \alpha/\varepsilon - A_{10}A_{01}/\varepsilon.
\label{eq:prod0d}
\end{eqnarray}
To verify (\ref{eq:sum0d}) 
we consider the expression $\alpha$ given
in (\ref{eq:alphadef}).
We simplify this expression by evaluating 
$A_{11}$ in terms of
 $\theta_0, \theta_d, A_{00}, A_{dd}$ and
$\theta^*_0, \theta^*_1, \ldots, \theta^*_d$ using
(\ref{eq:aiiform}) and Theorem 
\ref{lem:step2}.  Simplifying the result further using
(\ref{eq:epsdef2})
we find $\alpha = \varepsilon(\theta_0+\theta_d-A_{00})$ and line 
 (\ref{eq:sum0d}) follows. 
 To verify
(\ref{eq:prod0d}) 
we evaluate
 the product $A_{10}A_{01}$ in terms of
 $\theta_0, \theta_d, A_{00}, A_{dd}$ and
$\theta^*_0, \theta^*_1, \ldots, \theta^*_d$ 
using 
(\ref{eq:crossprod}) and 
Theorem \ref{lem:step2}.
Simplifying the result using
(\ref{eq:epsdef1})
 we obtain 
$A_{10}A_{01}=-\varepsilon (A_{00}-\theta_0)(A_{00}-\theta_d)$.
Combining this with 
(\ref{eq:sum0d}) we routinely obtain
(\ref{eq:prod0d}).
\hfill $\Box $\\

\section{Transition matrices and polynomials}

\medskip
\noindent 
Let $\Phi$ denote a Leonard system over $\K$ and let
$(\theta_i, \theta^*_i, i=0..d;  \varphi_j, \phi_j, j=1..d)$
denote the corresponding parameter array.
Let $\alg$
denote the ambient algebra of $\Phi$.
Let $\flat :
\alg \rightarrow 
 \hbox{Mat}_{d+1}(\K)$
 denote
 the  $TD$-$D$ canonical map for $\Phi$, from 
Definition
\ref{def:flatcon}.
Let $\sharp
: \alg \rightarrow 
 \hbox{Mat}_{d+1}(\K)$
 denote the 
 $TD$-$D$ canonical map for $\Phi^*$.
We describe how  
  $\flat $ and 
 $\sharp$ are related. 
To do this we cite some facts from
\cite[Section 16]{LS24}.
For $0 \leq i,j\leq d$  we define
the scalar
\begin{eqnarray}
{\mathcal P}_{ij}=
\sum_{n=0}^d
\frac{(\theta_i-\theta_0)(\theta_i-\theta_1)\cdots (\theta_i-\theta_{n-1})
(\theta^*_j-\theta^*_0)(\theta^*_j-\theta^*_1)
\cdots (\theta^*_j-\theta^*_{n-1})
}{\varphi_1 \varphi_2 \cdots \varphi_n}.
\label{eq:sumpart}
\end{eqnarray}
%We remark ${\mathcal P}_{ij}$ is a polynomial of degree $j$ in
%$\theta_i$ and a polynomial of degree $i$ in $\theta^*_j$.
Let $P$ denote the matrix in 
 $\hbox{Mat}_{d+1}(\K)$
which has entries
\beast
P_{ij}=k_j{\mathcal P}_{ij} \qquad \qquad  
(0 \leq i,j\leq d),
\eeast
where ${\mathcal P}_{ij}$ is from
(\ref{eq:sumpart})
and where
 $k_j$ equals
\beast
\frac{\varphi_1 \varphi_2 \cdots \varphi_j}{\phi_1 \phi_2 \cdots \phi_j}
\eeast
times
\beast
\frac{(\theta^*_0-\theta^*_1)(\theta^*_0-\theta^*_2)\cdots
(\theta^*_0-\theta^*_d)}
{(\theta^*_j-\theta^*_0)\cdots (\theta^*_j-\theta^*_{j-1})
(\theta^*_j-\theta^*_{j+1})\cdots (\theta^*_j-\theta^*_d)}
\eeast
for $0 \leq j \leq d$. Then $P_{i0}=1$ for $0 \leq i \leq d$ and
$X^\sharp P = PX^\flat $ for all $ X \in 
\alg $.
Let $P^*$ denote the matrix in
 $\hbox{Mat}_{d+1}(\K)$ which has entries
\beast
P^*_{ij}=k^*_j{\mathcal P}_{ji} \qquad \qquad  
(0 \leq i,j\leq d),
\eeast
where ${\mathcal P}_{ji}$ is from 
(\ref{eq:sumpart})
and 
$k^*_j$ equals
\beast
\frac{\varphi_1 \varphi_2 \cdots \varphi_j}{\phi_d \phi_{d-1} \cdots \phi_{d-j+1}}
\eeast
times
\beast
\frac{(\theta_0-\theta_1)(\theta_0-\theta_2)\cdots
(\theta_0-\theta_d)}
{(\theta_j-\theta_0)\cdots (\theta_j-\theta_{j-1})
(\theta_j-\theta_{j+1})\cdots (\theta_j-\theta_d)}
\eeast
for $0 \leq j \leq d$.
Then $P^*_{i0}=1$ for $0 \leq i \leq d$ and
$X^\flat P^* = P^*X^\sharp $ for
all
 $ X \in 
\alg $.
Moreover $PP^*=\nu I$
where
\beast
\nu = 
\frac{
(\theta_0-\theta_1)(\theta_0-\theta_2)\cdots
(\theta_0-\theta_d)
(\theta^*_0-\theta^*_1)(\theta^*_0-\theta^*_2)\cdots
(\theta^*_0-\theta^*_d)}
{\phi_1 \phi_2 \cdots \phi_d}.
\eeast

\noindent 
We comment on 
(\ref{eq:sumpart}).
For
$0 \leq i,j\leq d$,
${\mathcal P}_{ij}$ is a polynomial of degree $j$ in
$\theta_i$ and a polynomial of degree $i$ in $\theta^*_j$.
The class of polynomials which can be obtained from
a parameter array
in this fashion
coincides with the class of polynomials which are contained in
the Askey scheme \cite{KoeSwa} 
and which are orthogonal with respect
to a measure which has finitely many nonzero values.
This class consists of the Krawtchouk,
Hahn, dual Hahn, Racah, the
$q$-analogs of these, and some polynomials obtained from
the $q$-Racah by letting $q = -1$.
See \cite[Appendix A]{LS99}
and \cite[p. 260]{BanIto}
for more details. 
To illustrate this we obtain
some Krawtchouk and $q$-Racah polynomials
from 
 the parameter arrays
given in Example 
\ref{ex:pa1}
and Example
\ref{ex:pa2}, respectively.

\begin{example} \cite[Section 16]{LS24}
\label{ex:KrawP}
\rm
Let 
$(\theta_i, \theta^*_i, i=0..d;  \varphi_j, \phi_j, j=1..d)$
denote the parameter array in
Example
\ref{ex:pa1}.
 Referring to the discussion in the first part of this section,
 for $0 \leq i,j \leq d$ we have 
\begin{eqnarray}
 {\mathcal P}_{ij} =
\sum_{n=0}^d \frac{(-i)_n (-j)_n 2^n}{(-d)_n n! }
%\,\frac{2^n}{n!},
\label{eq:2F1expand}
\end{eqnarray}
 where
\beast
(a)_n:=a(a+1)(a+2)\cdots (a+n-1) \qquad \qquad n=0,1,2,\ldots
\eeast
Moreover
\beast
k_j = \Biggl({{ d }\atop {j}}\Biggr), \qquad \qquad 
k^*_j = \Biggl({{ d }\atop {j}}\Biggr) \qquad \qquad 
(0 \leq j \leq d)
\eeast
and $\nu=2^d $.
We have $P=P^*$ and $P^2=2^dI$. 
For $0 \leq i,j\leq d$ 
the expression on the right
in 
(\ref{eq:2F1expand}) is equal to the hypergeometric series
\begin{eqnarray}
{{}_2}F_1\Biggl({{-i, -j}\atop {-d}}\;\Bigg\vert \;2\Biggr).
\label{eq:2F1not}
\end{eqnarray}
From this we find ${\mathcal P}_{ij}$ 
is a Krawtchouk polynomial of degree $j$ in
$\theta_i$ and a Krawtchouk polynomial of degree  $i$ in
$\theta^*_j$.
\end{example}

\begin{example} \cite[Section 16]{LS24}
\label{ex:qraclasttime}
\rm
Let 
$(\theta_i, \theta^*_i, i=0..d;  \varphi_j, \phi_j, j=1..d)$
denote the parameter array in
Example
\ref{ex:pa2}.
 Referring to the discussion in the first part of this section,
 for $0 \leq i,j \leq d$ we have 
\begin{eqnarray}
 {\mathcal P}_{ij} =
\sum_{n=0}^d \frac{(q^{-i};q)_n (sq^{i+1};q)_n 
(q^{-j};q)_n (s^*q^{j+1};q)_n q^n}
{(r_1q;q)_n(r_2q;q)_n (q^{-d};q)_n(q;q)_n} 
%\,\frac{q^n}{(q;q)_n }, 
\label{eq:uihyper}
\end{eqnarray}
where 
%%%%%%%%%%%%%%%%%%%%%%%%%%%%
\beast
(a;q)_n := (1-a)(1-aq)(1-aq^2)\cdots (1-aq^{n-1})\qquad \qquad n=0,1,2\ldots 
\eeast
Moreover 
\beast
k_j &=& \frac{(r_1q;q)_j(r_2q;q)_j(q^{-d};q)_j(s^*q;q)_j(1-s^*q^{2j+1})}
{s^jq^j(q;q)_j(s^*q/r_1;q)_j(s^*q/r_2;q)_j(s^*q^{d+2};q)_j(1-s^*q)},
\\
k^*_j &=& \frac{(r_1q;q)_j(r_2q;q)_j(q^{-d};q)_j(sq;q)_j(1-sq^{2j+1})}
{s^{*j}q^j(q;q)_j(sq/r_1;q)_j(sq/r_2;q)_j(sq^{d+2};q)_j(1-sq)} 
\eeast
for $0 \leq j \leq d$ and
\beast
\nu = \frac{(sq^2;q)_d (s^*q^2;q)_d}{r^d_1q^d(sq/r_1;q)_d(s^*q/r_1;q)_d}. 
\eeast
For $0 \leq i,j\leq d$ the expression
on the right in
(\ref{eq:uihyper})
 is equal to the basic hypergeometric series
\beast 
 {}_4\phi_3 \Biggl({{q^{-i}, \;sq^{i+1},\;q^{-j},\;s^*q^{j+1}}\atop
{r_1q,\;\;r_2q,\;\;q^{-d}}}\;\Bigg\vert \; q,\;q\Biggr).
\eeast
By this 
 and since
$r_1r_2=s s^*q^{d+1}$ 
we find  ${\mathcal P}_{ij}$ 
is a  
$q$-Racah polynomial of degree $j$ in $\theta_i$ and
a $q$-Racah polynomial 
of degree $i$ in $\theta^*_j$.
\end{example}

\section{Directions for further research}

\noindent In this section we give some suggestions
for further research.

\begin{problem} Let $\Phi$ denote the Leonard system
from (\ref{eq:ourstartingpt}). Let $\alpha, \alpha^*, \beta,\beta^*$
denote scalars in $\K$ such that $\alpha\not=0$ and $\alpha^*\not=0$.
Recall the sequence
$(\alpha A+\beta I,
\alpha^* A^*+\beta^* I; E_i, E^*_i,i=0..d )
$
is a Leonard system in $\mathcal A$. 
In some cases this system is isomorphic to a relative of
$\Phi$;
 describe all the cases where this occurs.
\end{problem}

%\begin{problem}
%\label{prob:rtu}
%Let $d$ denote a nonnegative  integer and let
%$(\theta_i, \theta^*_i, i=0\ldots d;  \varphi_j, \phi_j, j=1\ldots d)$
%denote a parameter array over $\K$. 
%Let us assume the scalar $q$ from
%Definition .. is not 1 or -1. Then there exists a sequence of
%scalars
%$r,t,u,r^*,t^*,u^*,\rho$ taken from $\K$ such that
%\beast
%\theta_i &=& r + t \vartheta_i +  u \vartheta_{i+1},
%\\
%\theta^*_i &=& r^* + t^* \vartheta_i +  u^* \vartheta_{i+1} \qquad \qquad 
%\eeast
%for $0 \leq i \leq d$
%and
%\beast
%\varphi_i &=& \vartheta_i(\rho + t t^*\vartheta_{i-1} 
%+(t u^*+t^* u) \vartheta_i + u u^* \vartheta_{i+1}),
%\\
%\phi_i &=& \vartheta_i(\rho + t^* u\vartheta_{i-1} 
%+(t t^*+u u^*) \vartheta_i + t u^* \vartheta_{i+1})
%\eeast
%for $1 \leq i \leq d$. (The scalars $\vartheta_i $ are from
%(\ref{eq:vartheta}).)
%We remark  
%$\vartheta_{i-1}-\beta \vartheta_i + \vartheta_{i+1}$ is 
%independent of $i$ for $1 \leq i \leq d-1$.
%\end{problem}

\begin{problem} Let $d$ denote a nonnegative integer.
Find all Leonard pairs 
$A,A^*$ in $\hbox{Mat}_{d+1}(\K)$ which satisfy
the following two conditions:
(i) $A$ is 
irreducible tridiagonal;
(ii)  $A^*$ is lower bidiagonal with 
$A_{i,i-1}=1$ for $1 \leq i \leq d$.
\end{problem}

\begin{problem} Let $d$ denote a nonnegative integer.
Find all Leonard pairs $A,A^*$  in $\hbox{Mat}_{d+1}(\K)$
such that each of $A, A^*$ is
irreducible tridiagonal.
\end{problem}

\begin{problem} Let $d$ denote a nonnegative integer.
Find all Leonard pairs $A,A^*$ in  
$\hbox{Mat}_{d+1}(\K)$ which satisfy the following
two conditions:
(i) 
each of $A,A^*$ is irreducible tridiagonal;
(ii) there exists a diagonal matrix $H$ in 
$\hbox{Mat}_{d+1}(\K)$ such that
$A=HA^*H^{-1}$.
\end{problem}

\begin{problem}
\label{prob:spin}
Let $A,A^*$ denote the Leonard pair  from Definition
\ref{def:lpabs}.
Determine when does there 
exist invertible elements $U, U^*$ in 
$\mathcal A$ which satisfy (i)--(iii) below:
(i) $UA=AU$; (ii) $U^*A^*=A^*U^*$;
(iii) $UA^*U^{-1}= U^{*-1}AU^{*}$.
%(iv) $U^{-1}A^*U= U^{*}AU^{*-1}$.
This problem arises naturally in the context of
a spin model contained in a Bose-Mesner algebra of
$P$- and $Q$-polynomial type
\cite{Curspin}.
\end{problem}

\begin{problem}
Let $V$ denote a vector space over $\K$ with finite positive
dimension. By a Leonard triple on $V$, we mean a three-tuple
of linear transformations 
$A:V\rightarrow V$, $A^*:V\rightarrow V$, 
 $A^\varepsilon:V\rightarrow V$ 
which  
 satisfy conditions (i)--(iii) below. 
\begin{enumerate}
\item There exists a basis for $V$ with respect to which
the matrix representing $A$ is 
diagonal and the 
the matrices representing $A^*$ and $A^{\varepsilon}$ are each irreducible
tridiagonal.
\item There exists a basis for $V$ with respect to which
the matrix representing $A^*$ is 
diagonal and the 
the matrices representing $A$ and $A^{\varepsilon}$ are each irreducible
tridiagonal.
\item There exists a basis for $V$ with respect to which
the matrix representing $A^{\varepsilon}$ is 
diagonal and the 
the matrices representing $A$ and $A^*$ are each irreducible
tridiagonal.
\end{enumerate}
Find all the Leonard triples. 
\end{problem}

\begin{remark}
Referring to Problem 
\ref{prob:spin},
let $A^\varepsilon$
denote the common value of $UA^*U^{-1}$, $U^{*-1}AU^{*}$.
%or the common value of $U^{-1}A^*U$, $U^{*}AU^{*-1}$.
Then $A,A^*, A^\varepsilon$ is a Leonard triple.
\end{remark}

\begin{conjecture}
\label{conj:w1111}
Let $\Phi$ denote the Leonard system
from
(\ref{eq:ourstartingpt}) and let $I$ denote
the identity element of $\mathcal A$. Then
 for all
$X \in {\mathcal A}$ the following are equivalent:
(i) both
\beast
&&E_iX E_j = 0  \quad 
\hbox{if} \quad |i-j|>1, 
\qquad \qquad (0 \leq i,j\leq d), 
\\
&&E^*_iX E^*_j = 0 \quad  \hbox{if} \quad |i-j|>1,
\qquad \qquad (0 \leq i,j\leq d);
\eeast
(ii) $X$ is a $\K$-linear combination of 
$I, A, A^*, AA^*, A^*A.$
\end{conjecture}

\begin{conjecture}
Let $\Phi$ denote the Leonard system from  
(\ref{eq:ourstartingpt}). Then for $0 \leq  r \leq d$
the elements
\beast
E^*_0, E^*_1, \ldots, E^*_r, E_r, E_{r+1}, \ldots, E_d
\eeast
together generate $\mathcal A$.
%Moreover the following products
%together form a basis for $\cal A$.
%\beast
%&& E^*_iE_dE^*_0E_j \qquad \qquad (0  \leq i\leq r, \qquad  r+1\leq j \leq d),
%\\
%&& E^*_iE_dE^*_j \qquad \qquad (0  \leq i\leq r, \qquad 0\leq j \leq r),
%\\
%&& E_iE^*_0E_j \qquad \qquad (r+1  \leq i\leq d, \qquad r+1\leq j \leq d),
%\\
%&& E_iE^*_0E_dE^*_j \qquad \qquad (r+1  \leq i\leq d, \qquad 0\leq j \leq r).
%\eeast
\end{conjecture}

\section{Acknowledgement} The author thanks
John Caughman and Hjalmar Rosengren for a conversation which inspired
Problem 
\ref{prob:spin}.
The author thanks 
Brian Curtin, Mark MacLean, and Raimundas Vidunas
for giving this manuscript a close reading
and offering many valuable suggestions.

%\small
%\bibliographystyle{plain}
%
%\bibliography{master}
%\normalsize
%\medskip
%
%\noindent Paul Terwilliger, Department of Mathematics, University of
%Wisconsin, 480 Lincoln Drive, Madison, Wisconsin, 53706, USA \hfil\break
%email: terwilli@math.wisc.edu \hfil\break
%
%

\noindent Paul Terwilliger  \hfil \break
\noindent Department of Mathematics \hfil \break
\noindent University of
Wisconsin \hfil \break
\noindent 480 Lincoln Drive \hfil\break
\noindent Madison, Wisconsin,  53706 USA \hfil\break
Email: terwilli@math.wisc.edu \hfil\break

\end{document}